\documentclass[reqno,tbtags,11pt]{amsart}

\usepackage[margin=1.25in]{geometry}

\usepackage{microtype}
\usepackage{mathtools}
\usepackage{array}
\usepackage{csquotes}
\usepackage{graphicx}
\usepackage{caption}
\usepackage{subcaption}
\usepackage[section]{placeins}

\newcommand{\dd}{\mathop{}\!{\mathrm {d}}}
\newcommand{\ii}{\mathrm {i}}

\DeclarePairedDelimiter{\abs}{\lvert}{\rvert}
\DeclarePairedDelimiter{\measure}{\lvert}{\rvert}
\DeclarePairedDelimiter{\norm}{\lVert}{\rVert}
\DeclarePairedDelimiter{\ip}{\langle}{\rangle}
\DeclarePairedDelimiter{\set}{\lbrace}{\rbrace}

\DeclareMathOperator{\sech}{sech}
\newcommand{\R}{\mathbb {R}}

\usepackage[hidelinks]{hyperref}
\hypersetup{
 pdftitle={Finding geometric criteria for negative eigenvalues of the Neumann--Poincaré operator},
 pdfauthor={Magnus Goffeng, Johan Helsing, Mikael Sundqvist}}

\newtheorem{theorem}{Theorem}[section]
\newtheorem{lemma}[theorem]{Lemma}
\newtheorem{proposition}[theorem]{Proposition}
\newtheorem{corollary}[theorem]{Corollary}

\theoremstyle{remark}
\newtheorem{remark}[theorem]{Remark}
\newtheorem{example}[theorem]{Example}

\numberwithin{equation}{section}

\author{M. Goffeng}
\author{J. Helsing}
\author{M. Sundqvist}

\address{Department of Mathematics, Lund University, Sweden}

\email{magnus.goffeng@math.lth.se}
\email{johan.helsing@math.lth.se}
\email{mikael.persson\_sundqvist@math.lth.se}


\title[Negative eigenvalues of the NP operator]
  {Finding geometric criteria for negative eigenvalues\\ of the Neumann--Poincaré operator}

\begin{document}

\begin{abstract}
  Motivated by a question of Grigori Rozenblum, we study how the geometry of the
  boundary of a convex domain $\Omega\subset\R^3$ affects whether the
  Neumann--Poincaré operator on $\partial\Omega$ has negative eigenvalues. First,
  we show that a flat face on the boundary of a bounded, convex Lipschitz domain
  produces an infinite-dimensional negative spectral subspace. For convex domains
  with $C^{1,\alpha}$ boundary, we bound the number of negative eigenvalues from
  below by the Dirichlet eigenvalue counting function of a projected patch over
  which the chords are sufficiently short. A Rellich identity gives a sufficient
  condition on a defining function to exclude negative eigenvalues. We also
  report numerical computations, which led us to the theoretical results.
\end{abstract}

\maketitle

\section{Introduction}

The Neumann--Poincar{\'e} (NP) operator is a classical boundary integral operator
associated with the Laplace equation. It governs, among other things, the
solvability of transmission problems and the electrostatic response of an
inclusion~\cite{AmmariKang2007,Kang2014,KangKimLeeShinYu2016}; a general account
of its spectral theory is given
in~\cite{AmmariMillienRuizZhang2017,AKMP21,Grieser2014}. In this paper we study
how the geometry of the boundary of a convex domain controls the sign of the
eigenvalues of its NP operator in dimension three. For smooth boundaries, there
are infinitely many positive eigenvalues~\cite{MR3D}, and the importance of the
negative eigenvalues is discussed in~\cite{Kang2022}. By smooth we mean
$C^\infty$. The motivation for this paper is the following question posed to the
authors by Grigori Rozenblum: \emph{Can we determine from the geometry of a
sufficiently regular, strictly convex domain whether there are negative
eigenvalues of its NP operator?} The natural follow-up question is: \emph{If yes,
can we determine how many there are?}

We give partial answers to these questions in terms of three geometric criteria
on the bounded Lipschitz domain $\Omega\subset\R^3$. First, if there is a plane
$\Pi$ such that $\Pi\cap\partial\Omega$ has nonempty interior in $\partial\Omega$
(that is, a flat face) and $\Omega$ is convex, then the NP operator has an
infinite-dimensional negative spectral subspace.
Second, if the chords in some direction are short over a large planar region $D$,
then $N_{-}(K^*)$ is bounded below by a Dirichlet counting function on $D$.
Third, the balanced-Hessian condition $(\Delta F)I\geq 2D^2F$ on a
defining function $F$ guarantees that every eigenvalue is nonnegative.
These criteria go beyond the spheroids, whose spectra are explicit, but they
neither characterize when negative eigenvalues exist nor give sharp bounds.

The asymptotic behavior of the eigenvalues of the NP operator is well understood
for smooth $\partial\Omega$. For a smooth surface, Miyanishi and
Rozenblum~\cite{MR3D} obtained separate asymptotic formulas for the positive and
negative eigenvalues in terms of the principal curvatures and proved the Weyl law
for their combined counting function. In particular, if~$\partial\Omega$ has
everywhere positive principal curvatures (with the convention that the sphere has
positive curvatures), only finitely many negative eigenvalues can occur, so it
makes sense to ask how many there are. This curvature assumption is stronger than
strict convexity. As a partial converse, Miyanishi and Rozenblum~\cite{MR3D}
proved that if there is a point of nonconvexity on~$\partial\Omega$, then there
are infinitely many negative eigenvalues. For a different concavity criterion, Ji
and Kang proved the existence of a negative spectral value of the NP operator on
the boundary or on its spherical inversion under a suitable geometric
hypothesis~\cite{JiKang19}.

Ando, Ji, Kang, Kawagoe, Miyanishi, and Rozenblum gave a sufficient condition for
a surface with two flat faces to have a negative eigenvalue~\cite{AJKMKR}. They
applied it in particular to what we call the cheese below; see
Figure~\ref{fig:curves}(\subref{fig:cheesecurve}).

\subsection{Setup of the spectral problem}

Let $\Omega\subset\R^3$ be a bounded Lipschitz domain, with boundary
$\partial\Omega$. We let $n_x$ denote the outward unit normal at
$x\in\partial\Omega$. Further regularity of~$\partial\Omega$ is assumed where
needed and is stated in each result; for the numerical experiments the boundary
is at least $C^{1,1}$. We use the NP sign convention adopted in~\cite{AKMP21}
and normalize the spectrum to lie in the interval $[-1,1]$. The single-layer
operator is defined as
\[
  S\varphi(x)
 \coloneqq
 -\frac{1}{4\pi}
 \int_{\partial\Omega}\frac{\varphi(y)}{\abs{x-y}}\dd S_y.
\]
The normalized double-layer operator is defined by
\[
 K\varphi(x)
 \coloneqq
 \frac{1}{2\pi}
 \operatorname{p.v.}\!\int_{\partial\Omega}
 \frac{n_y\cdot(y-x)}{\abs{x-y}^3}\varphi(y)\dd S_y.
\]
The operators $S$ and $K$ are bounded on
$L^2(\partial\Omega)$. The adjoint of the double-layer operator is given by
\begin{equation}
K^*\varphi(x)
\coloneqq
 \frac{1}{2\pi}
 \operatorname{p.v.}\!\int_{\partial\Omega}
 \frac{n_x\cdot(x-y)}{\abs{x-y}^3}\varphi(y)\dd S_y .
 \label{eq:intro-np}
\end{equation}
In this normalization $K1=1$. We call $K^*$ the NP operator. The reader should
be aware that there is conflicting terminology and normalization in the literature,
as discussed below in Appendix~\ref{secwarning}. We write $S_{\partial\Omega}$,
$K _ {\partial\Omega}$, and $K_{\partial\Omega}^*$ for $S$, $K$, and $K^*$ when we
need to emphasize which surface we are working on.

The operators $S$ and $K$ satisfy the Plemelj symmetrization principle
\[
  KS=SK^*.
\]
In particular, the single-layer operator symmetrizes $K^*$ through the energy
form
\begin{equation}
 (\varphi,\psi)_*
 \coloneqq
 \ip{\varphi,-S\psi}.
 \label{eq:intro-energy-form}
\end{equation}
The operator $-S$ is positive in the form sense and injective with dense range on
$L^2(\partial \Omega)$, so we can consider the completion $\mathfrak{H}$ of
$L^2(\partial \Omega)$ with respect to the energy
form~\eqref{eq:intro-energy-form}. For a Lipschitz boundary, the single-layer
isomorphism identifies $\mathfrak{H}$ with $H^{-1/2}(\partial\Omega)$, with
equivalent norms; see~\cite[Chapters~6--8]{McLean2000}. We consider $K^*$ as an
operator on $\mathfrak{H}$. This energy-space formulation is the setting for the
variational spectral description in~\cite{KPS07}. For a Lipschitz boundary, $K^*$
is self-adjoint in the energy space $\mathfrak{H}$, and for a $C^{1,\alpha}$
boundary, $\alpha>0$, it is compact on both $\mathfrak{H}$ and
$L^2(\partial\Omega)$.

We note that we could equally well consider the operator $K$ as a self-adjoint
operator on the $L^2$-dual space of $\mathfrak{H}$, as it would have the same
spectrum as $K^*$. If $\partial\Omega$ is smooth, ellipticity of $S$ ensures not
only that the two compact operators $K^*$ and $K$ have the same spectrum, but
also that this spectrum is the same on the energy space $\mathfrak{H}$, on its
$L^2$-dual, and on $L^2(\partial\Omega)$.

For $x\not\in\partial\Omega$, we define the single-layer potential
\begin{equation}\label{singledoubleintro}
 \mathcal S\varphi(x)
 =
 -\frac{1}{4\pi}\int_{\partial\Omega}\frac{\varphi(y)}{\abs{x-y}}\dd S_y.
\end{equation}

\subsection{Main results}

The spectrum of the self-adjoint operator $K^*$ on
$\mathfrak{H}$ is known to be a subset of $[-1,1]$, with $1$ an eigenvalue. The
constant function is an eigenfunction of $K$, whereas $S^{-1}1$ is an
eigenfunction of $K^*$. For Lipschitz boundaries, we write $N_-(K^*)$ for the
dimension of the negative spectral subspace, that is,
\[
 N_-(K^*)\coloneqq\dim\operatorname{Ran}\mathbf 1_{(-\infty,0)}(K^*).
\]
For $C^{1,\alpha}$ boundaries,
compactness implies that $N_-(K^*)$ is the number of negative eigenvalues counted
with multiplicity.

\begin{figure}[htbp]
    \centering
    \begin{subfigure}{0.32\textwidth}
        \centering
        \includegraphics[height=2.5cm,page=1]{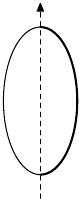}
        \caption{Prolate spheroid}
        \label{fig:prolatecurve}
    \end{subfigure}
    \hfill
    \begin{subfigure}{0.32\textwidth}
        \centering
        \includegraphics[height=2.5cm,page=2]{np-fig.pdf}
        \caption{Oblate spheroid}
        \label{fig:oblatecurve}
    \end{subfigure}
    \hfill
    \begin{subfigure}{0.32\textwidth}
        \centering
        \includegraphics[height=2.5cm,page=3]{np-fig.pdf}
        \caption{Superellipsoid}
        \label{fig:ellipsoidcurve}
    \end{subfigure}
    \\
    \begin{subfigure}{0.32\textwidth}
        \centering
        \includegraphics[height=2.5cm,page=4]{np-fig.pdf}
        \caption{Cheese}
        \label{fig:cheesecurve}
    \end{subfigure}
    \hfill
    \begin{subfigure}{0.32\textwidth}
        \centering
        \includegraphics[height=2.5cm,page=5]{np-fig.pdf}
        \caption{Sausage}
        \label{fig:sausagecurve}
    \end{subfigure}
    \hfill
    \begin{subfigure}{0.32\textwidth}
        \centering
        \includegraphics[height=2.5cm,page=6]{np-fig.pdf}
        \caption{Swedish coconut macaroon}
        \label{fig:coconutcurve}
    \end{subfigure}
    \caption{Generating curves (bold) and $C^{1,1}$ profiles (bold and thin) used
    to create surfaces of revolution. The dashed line is the axis of revolution.
    The spheroids are smooth calibration examples, the superellipsoid has
    curvature degeneracies, and the joined profiles are assembled from straight
    and curved parts. Negative eigenvalues were found for all examples
    except (a) and (e). See Section~\ref{sec:numerics} for more details.}
    \label{fig:curves}
\end{figure}

Before stating our main results, let us describe the ideas that led to them.
There are few explicitly computed examples, except for spheroids---the surfaces
obtained by rotating an ellipse about one of its axes, as in
Figures~\ref{fig:curves}(\subref{fig:prolatecurve})
and~\ref{fig:curves}(\subref{fig:oblatecurve}). For spheroids, one sees a
contrast between long prolate spheroids with no negative eigenvalues and oblate
spheroids, which have negative eigenvalues when they are thin enough. The
distance from a sphere alone does not decide the sign; the type of deformation
also matters.

To further explore what geometric features produce negative
eigenvalues, we carried out several numerical experiments. Their outcome is
described in more detail below in Section~\ref{sec:numerics}. The numerical
experiments concerned surfaces of revolution generated by the curves in
Figure~\ref{fig:curves}. Computations on the prolate spheroid and sausage reveal
no negative eigenvalues, whereas for the oblate spheroid and superellipsoid some
eigenvalue branches (eigenvalues as functions of a shape parameter) become
negative. The three joined profiles are
$C^{1,1}$: the cheese profile consists of two semicircles joined by straight
segments, the sausage is obtained by rotating the cheese profile about a
different axis, and the Swedish coconut macaroon profile is a rounded triangle.
The cheese and the Swedish coconut macaroon have flat faces, while
the sausage has no flat face.

From the experiments, we formed an intuition that negative eigenvalues arise when
$\partial\Omega$ has close-to-flat faces. There are likely to
be further geometric features of relevance: wedges~\cite{Perfekt19} can produce
negative essential spectrum, while conical points produce essential spectral
intervals whose location depends on the cone geometry~\cite{HelsingPerfekt18}.
Already for smooth $\partial\Omega$, thin oblate spheroids acquire many negative
eigenvalues as they flatten and develop a wedge-like feature, whereas the spectra
of increasingly elongated prolate spheroids contain no negative eigenvalues and
become densely distributed in $[0,1]$ in our normalization~\cite{AKLM22}. Our
first main result supports the intuition that close-to-flat faces contribute to
the negative spectrum:

\begin{theorem}\label{thm:flat}
  Let $\Omega\subset\R^3$ be a bounded, convex domain, and suppose that
  $\partial\Omega$ contains a flat face, i.e., there is a plane $\Pi$ such that
  $\Pi\cap\partial\Omega$ has nonempty interior in $\partial\Omega$. Then
  \[
    N_-(K^*)=\infty.
  \]
  In particular, if $\partial\Omega$ is additionally $C^{1,\alpha}$, then $K^*$ has
  infinitely many negative eigenvalues, counted with multiplicity.
\end{theorem}

The reader can find the proof of Theorem~\ref{thm:flat} in
Section~\ref{sec:flat}. A natural question is the following: if
$\Omega\subset\R^3$ is a smooth, bounded, convex domain whose boundary contains a
flat face and has positive principal curvatures away from it, what are the
spectral asymptotics of the negative eigenvalues of $K^*$?

\begin{figure}[tbp]
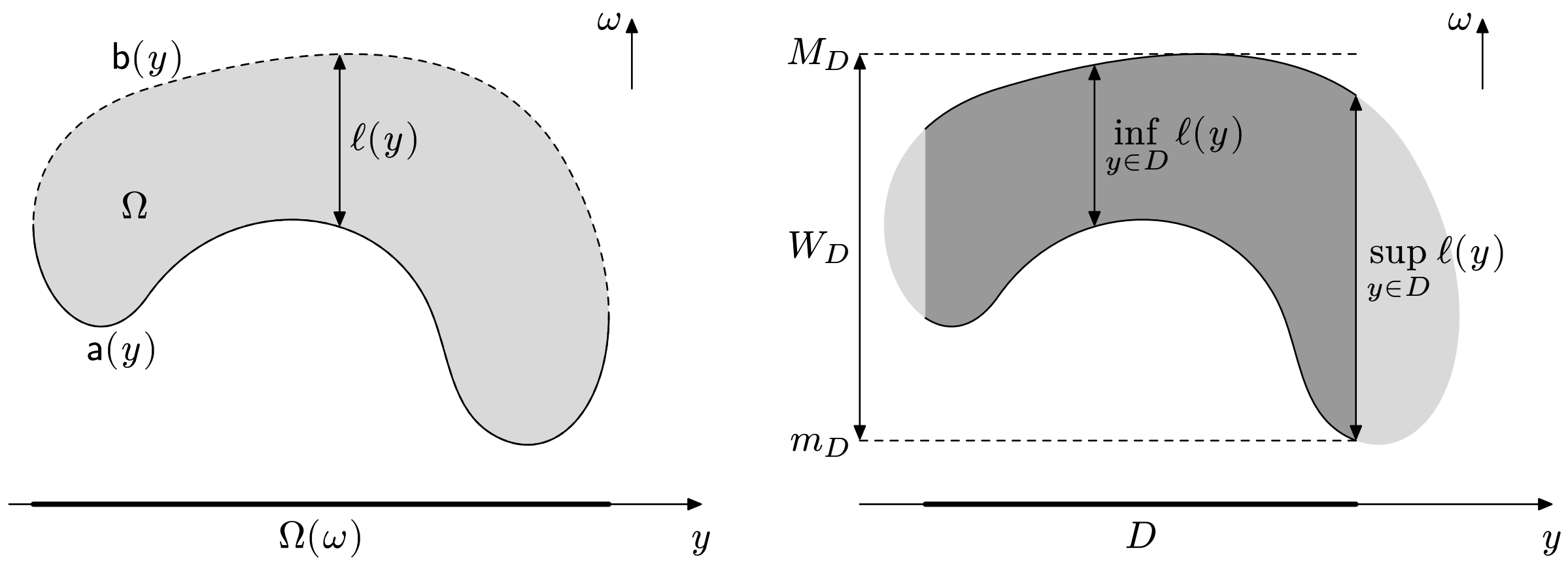

    \centering
    \begin{subfigure}{0.48\textwidth}
        \centering
        \includegraphics[page=7]{np-fig.pdf}
        \label{fig:Omegaa}
    \end{subfigure}
    \hfill
    \begin{subfigure}{0.48\textwidth}
        \centering
        \includegraphics[page=8]{np-fig.pdf}
        \label{fig:Omegab}
    \end{subfigure}
    \caption{Schematic (two-dimensional) figures showing the projected-chord
    geometry and the quantities defined
    in~\eqref{geomsetup1}--\eqref{geomsetup2}.}
    \label{fig:Omegas}
\end{figure}

We can also give a partial answer to Rozenblum's question. We require some
geometric setup for $\Omega$, and refer also to the sketches in
Figure~\ref{fig:Omegas} for the notation that is about to be introduced. Fix a
unit vector $\omega\in \R^3$. We write
\[
  \Omega(\omega)\coloneqq\operatorname{proj}_{\omega^\perp}\Omega
\]
for the projection of $\Omega$ onto the plane $\omega^\perp$. We say that
$\Omega$ \emph{satisfies the chord condition along~$\omega$} if there are
functions $\mathsf{a},\mathsf{b}\colon\Omega(\omega)\to \R$ such that
\[
 \Omega=\set{y+t\omega:y\in \Omega(\omega),\, \mathsf{a}(y)<t<\mathsf{b}(y)}.
\]
We note that $\mathsf{a}$ and $\mathsf{b}$ are uniquely determined by this
property. If $\Omega$ is a convex domain, then $\Omega$ satisfies the chord
condition along any direction $\omega$. We call a nonempty, bounded Lipschitz
domain $D\subset\Omega(\omega)$ with $\overline D\subset\Omega(\omega)$ a
\emph{projected patch}. We introduce the notation
\begin{equation}
\label{geomsetup1}
  \ell(y)=\mathsf{b}(y)-\mathsf{a}(y),\quad
  M_D=\sup_{y\in D} \mathsf{b}(y),\quad
  m_D=\inf_{y\in D} \mathsf{a}(y),\quad
  W_D=M_D-m_D.
\end{equation}
For each projected patch $D$, we also assume that $\mathsf{a}$ and $\mathsf{b}$
are Lipschitz on a neighborhood of $\overline D$. We will state a result in terms
of the quantities
\begin{equation}
\label{geomsetup2}
 h_D=\sup_{y\in D}\ell(y),
 \quad\text{and}\quad
 g_D = W_D - \inf_{y\in D} \ell(y).
\end{equation}
Here $h_D$ is the supremum of the chord lengths over $D$, while $g_D$ is the
difference between the width $W_D$ of $\Omega$ over $D$ and the infimum of the
chord lengths over~$D$.

We let $\mu_k(D)$ denote the $k$th Dirichlet eigenvalue of the positive planar
Laplacian on $D$, counted with multiplicity, and define the counting function
\begin{equation}
\label{defknefknd}
 N_D^{\mathrm{Dir}}(\Lambda)
 \coloneqq\#\set{j\in\mathbb N:\mu_j(D)<\Lambda}.
\end{equation}
Our second main result bounds $N_-(K^*)$ from below in terms of the counting
function $N_D^{\mathrm{Dir}}$.

\begin{theorem}\label{thm:counting}
Let $\Omega\subset\R^3$ be a bounded Lipschitz domain and $\omega\in \R^3$ a unit
vector along which $\Omega$ satisfies the chord condition. Let $D\subset
\Omega(\omega)$ be a bounded Lipschitz domain with $\overline D\subset
\Omega(\omega)$ such that $\mathsf{a}$ and $\mathsf{b}$ are Lipschitz near
$\overline D$. Then we have the lower bound
\[
  N_-(K^*)\geq N_D^{\mathrm{Dir}}\left(\frac{16}{(h_D+\sqrt{h_D}\sqrt{h_D+4g_D})^2}\right).
\]
\end{theorem}

Theorem~\ref{thm:counting} follows from Lemma~\ref{ljnkajdnakjnd} in
Section~\ref{sec:chords} and makes quantitative what we observed numerically for
thin oblate spheroids. Recall that a spheroid with semi-axes $(a,a,b)$ is defined
as the set of points $x=(x_1,x_2,x_3)\in\R^3$ satisfying the equation
\begin{equation}\label{eq:spheroidequation}
 \frac{x_1^2+x_2^2}{a^2}+\frac{x_3^2}{b^2}=1.
\end{equation}
It is called oblate if $a>b$ and prolate if $b>a$.

In Corollary~\ref{cor:flattening-count} we give a bound where we instead use the
width of $\Omega$ in the direction $\omega$, and in
Corollary~\ref{cor:oblate-flattening-count} we apply it to oblate spheroids and
find an asymptotic lower bound on $ N_-(K^*)$ as the aspect ratio $R = a/b \to
\infty$.

Our third main result gives a sufficient condition for the absence of negative
eigenvalues.

\begin{theorem}
\label{thm:balanced-hessian}
Suppose that the bounded domain $\Omega$ is given by
\[
 \Omega=\set{x\in\R^3:F(x)<0},
\]
where $F\in C^2(\R^3)$ has $0$ as a regular value and $D^2F$ is bounded. If
$(\Delta F)(x)I-2D^2F(x)\geq0$ as a quadratic form for every $x\in\R^3$, then
$K^*$ has no negative eigenvalues.
\end{theorem}

\begin{remark}
\label{alkknaldnad}
For an ellipsoid the condition $(\Delta F)I-2D^2F\geq 0$ becomes the triangle
inequality for the reciprocal squares of the semi-axes. It is automatic for
prolate spheroids, and Theorem~\ref{thm:balanced-hessian} gives the sufficient
condition $a/b\leq\sqrt{2}$ for an oblate spheroid with semi-axes $(a,a,b)$ to
have no negative eigenvalues. In the opposite direction, the projected-chord
criterion in Example~\ref{ex:oblate-chord} gives a negative eigenvalue for an
oblate spheroid once $a/b>7/\sqrt{6}$. Neither of the thresholds
$\sqrt2\approx1.41$ and $7/\sqrt6\approx2.86$ is sharp.
\end{remark}

\subsection{The axial dipole threshold for oblate spheroids}

The results above explain many of the phenomena discovered numerically and
discussed in Section~\ref{sec:numerics}. However, they do not answer Rozenblum's
question completely. We summarize here the axial dipole threshold for oblate
spheroids.

The threshold is based on exact separation-of-variables formulas for spheroids
due to Ahner and Arenstorf~\cite{AhnerArenstorf86}. Ahner later gave a detailed
treatment of the oblate case~\cite{Ahner}. Ahner's calculation exhibited a
negative eigenvalue at finite aspect ratio; the authors of~\cite{AKLM22} describe
this as the first known three-dimensional example and note that the negative sign
was found numerically. For an oblate spheroid with aspect ratio $R=a/b>1$, the
axial dipole branch $(m,n)=(0,1)$ (with azimuthal index $m$ and degree $n$; see
Section~\ref{sec:oblate-exact}) has a unique zero $R_*>1$. The numerical
comparisons in Section~\ref{sec:oblate-exact} suggest that this is the first zero
crossing among all branches; we do not prove that ordering here.

The axial dipole eigenvalue is
\[
 \lambda_{0,1}(R)=1-2L_3(R),
 \qquad
 L_3(R)=\frac1{e^2}\biggl(1-
 \frac{\sqrt{1-e^2}}{e}\arcsin e\biggr),
 \qquad e=\sqrt{1-R^{-2}},
\]
where $L_3$ is the depolarization factor in the short-axis
direction~\cite{DiFratta2016}. A direct calculation shows that $L_3$ is strictly
increasing in $R$ from $1/3$ (as $R\downarrow1$) to $1$ (as $R\to\infty$), so the
unique zero of this branch is characterized by
\begin{equation}\label{eq:Rstar}
 L_3(R_*)=\frac12,
 \qquad
 R_*=1.815270959479\ldots\,.
\end{equation}

In particular, negative eigenvalues exist whenever $R>R_*$. This dipole
calculation alone does not exclude negative eigenvalues on other branches for
$1<R\leq R_*$. We have reason to suspect that $K^*$ has no negative eigenvalues
in this range. ChatGPT supplied a rather technical argument for this claim, which
we do not include here. A comparable characterization of negative spectrum for
general smooth convex domains remains an open question.

\subsection{Contents of the paper}

Section~\ref{sec:preliminaries} develops the energy test. It expresses the NP
quadratic form through the difference between the interior and exterior Dirichlet
energies of common boundary data, and then applies the min--max principle to
estimate the dimension of the negative spectral subspace.

Section~\ref{sec:flat} applies this test to domains with a flat face. A
comparison with half-space harmonic extensions shows that a flat face produces
an infinite-dimensional negative spectral subspace.

Section~\ref{sec:chords} introduces the projected-chord criterion and proves a
lower bound on $N_{-}(K^{*})$ in terms of a Dirichlet eigenvalue counting
function. The section also gives explicit consequences for surfaces of
revolution, including the flattening bound for oblate spheroids.

Section~\ref{sec:positivity} derives a Rellich identity and the balanced-Hessian
condition, which guarantees nonnegativity of the NP spectrum. Its applications
to ellipsoids and spheroids are followed by a stability result for perturbations
of the defining function.

Section~\ref{sec:numerics} discusses numerical eigenvalue computations for the
examples in Figure~\ref{fig:curves}, including a comparison with results from
Ahner's exact formulas for oblate spheroids.

Appendix~\ref{secwarning} records the sign and normalization conventions for the
single- and double-layer operators. This is useful when comparing the spectral
statements below with results in the literature that use a different convention.

\subsection*{Acknowledgments}

We thank Grigori Rozenblum for introducing us to the problem and for sharing the
inspiring manuscript~\cite{AJKMKR} with us.
We also thank Heiko Gimperlein for discussions.
This work was supported by the Swedish Research
Council under contracts VR~2025-03923 and VR~2021-03720.

\subsection*{Use of AI}

The work in this paper started from numerical experiments carried out by the
authors with the aim of determining the geometric features that produce negative
eigenvalues for the NP operator. These numerical experiments led us to form an
intuition for the problem. Building on this intuition, and with assistance from
OpenAI’s ChatGPT (version 5.6), we obtained the theoretical results. ChatGPT was
also used for drafting.

\section{The energy test}
\label{sec:preliminaries}

Our main tool for detecting negative spectrum is a well-known energy test, which
compares the quadratic form of~$K^*$ with interior and exterior Dirichlet
energies; see~\cite{KPS07,AKMP21}. Throughout this section, $\Omega$ denotes a
bounded Lipschitz domain in $\R^3$. We let $\mathcal
H^1(\R^3\setminus\overline\Omega )$ denote the space of functions $u$ such that
$u\in H^1((\R^3\setminus\overline\Omega) \cap B_R)$ for every $R>0$ and
\[
 \int_{\R^3\setminus\overline\Omega }
 \left(\abs{\nabla u}^2+\frac{\abs u^2}{1+\abs x^2}\right)\dd x<\infty.
\]
This weighted finite-energy space consists of functions that are $H^1$ up to
$\partial\Omega$. The exterior Dirichlet problem and its decay condition are
discussed in~\cite[Chapter~8, in particular Theorem~8.10]{McLean2000}; see
also~\cite[Section~1.6]{ChandlerWildeSpence2024}.

\begin{lemma}
\label{lem:miniinad}
Assume that $f\in H^{1/2}(\partial\Omega)$ and consider the affine spaces
\begin{align*}
  V^{\mathrm{in}}(f)&\coloneqq\set{u\in H^1(\Omega): \; u|_{\partial\Omega}=f} \quad\text{and}\\
  V^{\mathrm{out}}(f)&\coloneqq\set{u\in\mathcal H^1(\R^3\setminus\overline\Omega ): \; u|_{\partial\Omega}=f}.
\end{align*}
The energy functional $V^{\mathrm{in}}(f)\to \R$,
\[
  u\mapsto \int_\Omega\abs{\nabla u}^2\dd x,
\]
has a unique minimizer $u_{\mathrm{in}}$, which is the harmonic extension of $f$
to $\Omega$. The energy functional $V^{\mathrm{out}}(f)\to \R$,
\[
  u \mapsto
  \int_{\R^3\setminus\overline\Omega }\abs{\nabla u}^2\dd x,
\]
has a unique minimizer $u_{\mathrm{out}}$, which is the harmonic extension of
$f$ to $\R^3\setminus\overline\Omega $ satisfying the decay conditions
\[
 u_{\mathrm{out}}(x)=O(\abs{x}^{-1}),\qquad
 \nabla u_{\mathrm{out}}(x)=O(\abs{x}^{-2})\qquad(\abs{x}\to\infty).
\]
\end{lemma}

\begin{proof}
The single-layer isomorphism gives $\varphi=S^{-1}f$. The restrictions of
$u=\mathcal S\varphi$ to $\Omega$ and $\R^3\setminus\overline\Omega $ are
harmonic, have trace $f$, and belong to the stated spaces; the exterior
restriction has the stated decay~\cite[Chapters~6--8]{McLean2000}. For a
zero-trace variation $w$, weak harmonicity gives, with $G=\Omega$ or
$G=\R^3\setminus\overline\Omega $,
\[
 \int_G\abs{\nabla(u+w)}^2\dd x
 =\int_G\abs{\nabla u}^2\dd x+\int_G\abs{\nabla w}^2\dd x.
\]
In the exterior, the identity follows first for compactly supported variations
and then by a cutoff at infinity in $\mathcal H^1(\R^3\setminus\overline\Omega
)$. The identity proves minimality, and the equality case together with the zero
trace proves uniqueness.
\end{proof}

We use the following notation for the Dirichlet energies:
\begin{align*}
  E_{\mathrm{in}}(f)
  & \coloneqq \int_\Omega\abs{\nabla u_{\mathrm{in}}}^2\dd x
  =
  \inf_{u\in V^{\mathrm{in}}(f)}\int_\Omega\abs{\nabla u}^2\dd x,
  \\
  E_{\mathrm{out}}(f)
  &\coloneqq \int_{\R^3\setminus\overline\Omega}
  \abs{\nabla u_{\mathrm{out}}}^2\dd x
  =
  \inf_{u\in V^{\mathrm{out}}(f)}\int_{\R^3\setminus\overline\Omega}\abs{\nabla u}^2\dd x.
\end{align*}
We call an element $u\in V^{\mathrm{in}}(f)$ an \emph{admissible competitor for the
Dirichlet problem} in $\Omega$ and an element $u\in V^{\mathrm{out}}(f)$ an admissible
competitor for the Dirichlet problem in $\R^3\setminus\overline\Omega$. The following
energy test relates these energies to the negative spectrum of $K^*$.

\begin{lemma}
\label{lem:energy-test}
If there is a nonzero $f\in H^{1/2}(\partial\Omega)$ such that
\begin{equation}
 E_{\mathrm{in}}(f)>E_{\mathrm{out}}(f),
 \label{eq:energy-test}
\end{equation}
then $N_-(K^*)\geq1$. More generally, if the
inequality~\eqref{eq:energy-test} holds for every nonzero element of a
$k$-dimensional subspace of $H^{1/2}(\partial \Omega)$, then $N_-(K^*)\geq k$. If
$\partial\Omega$ is $C^{1,\alpha}$, these conclusions give negative eigenvalues
counted with multiplicity.
\end{lemma}

\begin{proof}
The single-layer operator is an isomorphism
\[
 S\colon H^{-1/2}(\partial\Omega)\to H^{1/2}(\partial\Omega).
\]
Write $f=S\varphi=(\mathcal{S}\varphi)|_{\partial\Omega}$, where $\mathcal{S}$ is
as in~\eqref{singledoubleintro}. In particular, we have by uniqueness of
minimizers in Lemma~\ref{lem:miniinad} that
\[
  u_{\mathrm{in}}
  =
  (\mathcal{S}\varphi)|_{\Omega}
  \quad\text{and}\quad
  u_{\mathrm{out}}=(\mathcal{S}\varphi)|_{\R^3\setminus\overline\Omega}.
\]
With both normal derivatives taken in the direction of the outward normal of
$\Omega$, the jump relations are
\begin{equation}\label{eq:jump-in-out}
 \partial_n(\mathcal{S}\varphi)\big|_{\mathrm{in}}
 =\frac12(K^*-I)\varphi,
 \qquad
 \partial_n(\mathcal{S}\varphi)\big|_{\mathrm{out}}
 =\frac12(K^*+I)\varphi.
\end{equation}
Green's formula gives
\[
 E_{\mathrm{in}}(f)
 =\frac12\norm{\varphi}_*^2-\frac12(K^*\varphi,\varphi)_*,
 \qquad
 E_{\mathrm{out}}(f)
 =\frac12\norm{\varphi}_*^2+\frac12(K^*\varphi,\varphi)_*.
\]
Consequently,
\begin{equation}\label{eq:energytest}
 (K^*\varphi,\varphi)_*
 =E_{\mathrm{out}}(f)-E_{\mathrm{in}}(f).
\end{equation}
The strict inequality in~\eqref{eq:energy-test} therefore gives a negative
Rayleigh quotient. The spectral theorem then gives $N_-(K^*)\geq1$ and, more
generally, $N_-(K^*)\geq k$ when the quadratic form is negative on a
$k$-dimensional subspace. In the compact case, this is the usual
eigenvalue-counting conclusion of the min--max principle.
\end{proof}

\begin{remark}
\label{remsummada}
For an eigenmode, the energy test is exact. If
$K^*\varphi=\lambda\varphi$ and $f=S\varphi$, then
\[
 E_{\mathrm{in}}(f)=\frac{1-\lambda}{2}\norm{\varphi}_*^2,
 \qquad
 E_{\mathrm{out}}(f)=\frac{1+\lambda}{2}\norm{\varphi}_*^2.
\]
Consequently,
\[
\lambda
 =\frac{E_{\mathrm{out}}(f)-E_{\mathrm{in}}(f)}{E_{\mathrm{in}}(f)+E_{\mathrm{out}}(f)}.
\]
Thus $\lambda<0$ exactly when $\mathcal{S}\varphi$ has more Dirichlet energy in
$\Omega$ than in $\R^3\setminus\overline\Omega$.

By Lemma~\ref{lem:miniinad}, any admissible competitor bounds
$E_{\mathrm {out}}$ from above. In Section~\ref{sec:chords} the data are
$\pm\eta$ above and below a projected patch, and
$E_{\mathrm {in}} \geq 4\int\eta^2/\ell$ is large when the chords are short.
\end{remark}

\section{Domains with a flat face}
\label{sec:flat}

As a first application of the energy test of the previous section,
summarized in Remark~\ref{remsummada}, we treat convex domains whose
boundary has a flat face.

\begin{lemma}
\label{prop:flat-face}
Let $\Omega\subset\R^3$ be a bounded, convex domain, and suppose that
$\Gamma_0\subset\partial\Omega$ is a flat face. Then there is an
infinite-dimensional subspace $\mathcal X\subset H^{1/2}(\partial\Omega)$,
consisting of functions supported in $\Gamma_0$, such that
\begin{equation}
 E_{\mathrm{in}}(f)>E_{\mathrm{out}}(f)
 \qquad\text{for every }0\ne f\in\mathcal X.
 \label{eq:flat-face-energy}
\end{equation}
Consequently, $N_-(K^*)=\infty$.
\end{lemma}

\begin{proof}
Let $\Pi$ be the plane containing $\Gamma_0$, and choose the half-spaces $H_-$
and $H_+$ so that $\Omega\subset H_-$ and $H_+$ lies on the exterior side of the
face. Let $\mathcal X$ be the space of extensions by zero to $\partial\Omega$ of
functions in $C_c^\infty(\Gamma_0)$. For $f\in\mathcal X$, regard $f$ as its zero
extension to $\Pi$, and let $u_H^-$ and $u_H^+$ be its harmonic extensions to
$H_-$ and $H_+$, respectively, given by the Poisson kernel and satisfying the
decay condition $u_H^\pm(x)=O(\abs{x}^{-1})$ as $\abs{x}\to\infty$. These
extensions minimize the Dirichlet energy in the corresponding weighted
finite-energy spaces with trace $f$, as follows from a Fourier transformation
parallel to $\Pi$. By reflection across $\Pi$, their Dirichlet energies are
equal; denote the common value by $E_H(f)$.

Extending $u_{\mathrm{in}}$ by zero from $\Omega$ to $H_-$ gives an admissible
competitor for the half-space Dirichlet problem, since its trace vanishes on
$\partial\Omega\cap H_-$. Hence $E_{\mathrm{in}}(f)\geq
E_H(f)$ by the half-space Dirichlet principle. For $f\ne0$, the inequality is
strict: equality would make this zero extension equal to $u_H^-$, but the former
vanishes on the nonempty open set $H_-\setminus\overline\Omega$. By unique
continuation, this would force $u_H^-$ to vanish identically, contrary to its
nonzero boundary trace~$f$. On the other hand, $H_+$ is contained in the exterior
of $\Omega$. Since $u_H^+$ has zero trace on $\Pi\setminus\overline\Omega$,
extending it by zero to the remainder of the exterior gives an admissible
competitor, and therefore $E_{\mathrm{out}}(f)\leq E_H(f)$ by
the exterior part of Lemma~\ref{lem:miniinad}. Combining the two inequalities
proves~\eqref{eq:flat-face-energy}.

The space $S^{-1}\mathcal X$ is
infinite-dimensional, and the energy test~\eqref{eq:energytest} gives
\[
 (K^*\varphi,\varphi)_*=E_{\mathrm{out}}(S\varphi)
 -E_{\mathrm{in}}(S\varphi)<0,
 \qquad 0\ne\varphi\in S^{-1}\mathcal X.
\]
Thus the quadratic form of the bounded self-adjoint operator $K^*$ on
$\mathfrak{H}$ has infinite negative index. If the range of its negative
spectral projection $\mathbf 1_{(-\infty,0)}(K^*)$ were finite-dimensional, the
infinite-dimensional space $S^{-1}\mathcal X$ would contain a nonzero $\varphi$
in its orthogonal complement, where the quadratic form is nonnegative; this
contradicts $(K^*\varphi,\varphi)_*<0$. Hence $N_-(K^*)=\infty$.
\end{proof}

\begin{proof}[Proof of Theorem~\ref{thm:flat}]
Lemma~\ref{prop:flat-face} gives $N_-(K^*)=\infty$ in the presence of a flat face.
If the boundary is $C^{1,\alpha}$, compactness
gives infinitely many negative eigenvalues counted with multiplicity.
\end{proof}

\section{Projected chords and negative eigenvalues}
\label{sec:chords}

In this section we use the chord notation introduced
around~\eqref{geomsetup1}--\eqref{geomsetup2} and Theorem~\ref{thm:counting}. The
estimates below use only this chord representation, not convexity in any other
direction. Thus Proposition~\ref{prop:projected-chord} and its counting argument also
apply when every line parallel to $\omega$ that meets $\Omega$ intersects it in
an interval, provided the endpoint functions are Lipschitz near $\overline D$ as
assumed above. With $D$, $\ell$, $M_D$, $m_D$, $W_D$, $h_D$, and $g_D$ as
in~\eqref{geomsetup1}--\eqref{geomsetup2} above, the defect
\[
 W_D-\ell(y)=\bigl(M_D-\mathsf{b}(y)\bigr)+\bigl(\mathsf{a}(y)-m_D\bigr)
\]
measures the variation of the upper and lower parts of $\partial\Omega$ over $D$. It vanishes when
the two parts lie in planes orthogonal to $\omega$ at constant separation. We
take real-valued test functions throughout this section. For $\eta\in H_0^1(D)$,
we use the notation
\[
\mathfrak{E}_D(\eta)\coloneqq\int_D\abs{\nabla\eta}^2\dd y.
\]

\begin{proposition}
\label{prop:projected-chord}
Assume that $\Omega$ is a bounded Lipschitz domain satisfying the chord condition
along $\omega$, and $D\subset \Omega(\omega)$ is a bounded Lipschitz domain with
$\overline D\subset \Omega(\omega)$. Assume that $\mathsf{a}$ and $\mathsf{b}$
are Lipschitz near $\overline D$. If there is a subspace $V\subset H^1_0(D)$
such that for every $0\ne\eta\in V$,
\begin{equation}
 4\int_D\frac{\eta(y)^2}{\ell(y)}\dd y
 >
 \int_D\bigl(W_D-\ell(y)\bigr)\abs{\nabla\eta(y)}^2\dd y
 +2\sqrt{\mathfrak{E}_D(\eta)}\norm{\eta}_{L^2(D)},
 \label{eq:projected-chord-condition}
\end{equation}
then $N_-(K^*)\geq\dim(V)$.
\end{proposition}

\begin{proof}
We again apply the energy test of Section~\ref{sec:preliminaries}.
For $0\ne\eta\in H^1_0(D)$, we prescribe the boundary data $f_\eta\in
H^{1/2}(\partial \Omega)$ defined by
\[
 f_\eta(y+\mathsf{b}(y)\omega)=\eta(y),
 \qquad
 f_\eta(y+\mathsf{a}(y)\omega)=-\eta(y),
\]
for $y\in D$, and set $f_\eta=0$ on the remainder of the boundary. The local
Lipschitz regularity of the upper and lower parts of $\partial\Omega$ over
$D$ and the zero trace of $\eta$ on $\partial D$ ensure that $f_\eta\in
H^{1/2}(\partial \Omega)$.

We estimate the two energies associated with $f_\eta$. Along a chord above $y\in
D$, the harmonic extension $u_{\mathrm{in}}$ of $f_\eta$ changes from
$-\eta(y)$ to $+\eta(y)$. The one-dimensional Cauchy--Schwarz inequality gives
\[
 \int_{\mathsf{a}(y)}^{\mathsf{b}(y)}
 \abs{\partial_\omega u_{\mathrm{in}}(y+t\omega)}^2\dd t
 \geq \frac{4\eta(y)^2}{\ell(y)}.
\]
After integration over $D$,
\[
 E_{\mathrm{in}}(f_\eta)
 \geq 4\int_D\frac{\eta^2}{\ell}\dd y.
\]

For the exterior energy, use an admissible competitor. Take an $L> 0$ to be
specified later. Above the upper part of $\partial\Omega$ over $D$, define
\[
 v_L(y,t)=
 \begin{cases}
 \eta(y),&\mathsf{b}(y)\leq t\leq M_D,\\
 \eta(y)\exp(-(t-M_D)/L),&t\geq M_D,
 \end{cases}
\]
and below the lower part of $\partial\Omega$ over $D$, define
\[
 v_L(y,t)=
 \begin{cases}
 -\eta(y),&m_D\leq t\leq \mathsf{a}(y),\\
 -\eta(y)\exp(-(m_D-t)/L),&t\leq m_D.
 \end{cases}
\]
Set $v_L=0$ elsewhere in the exterior. The constant and exponential pieces have
matching traces at the corresponding interfaces $t=m_{D}$ and $t=M_{D}$ and zero
lateral trace because $\eta\in H_0^1(D)$. Together with the exponential decay,
this gives $v_L\in V^-(f_\eta)$. The energy in the two vertical columns is
\[
 \int_D\bigl(W_D-\ell(y)\bigr)\abs{\nabla\eta(y)}^2\dd y,
\]
and the two exponential tails contribute
$L\mathfrak{E}_D(\eta)+\norm{\eta}_{L^2(D)}^2/L$. Since $u_{\mathrm{out}}$
minimizes the Dirichlet energy,
\[
 E_{\mathrm{out}}(f_\eta)
 \leq
 \int_D\bigl(W_D-\ell\bigr)\abs{\nabla\eta}^2\dd y
 +L\mathfrak{E}_D(\eta)+\frac{\norm{\eta}_{L^2(D)}^2}{L}.
\]
Choosing $L=\norm{\eta}_{L^2(D)}/\sqrt{\mathfrak{E}_D(\eta)}$ minimizes the last two
terms and gives $2\sqrt{\mathfrak{E}_D(\eta)}\norm{\eta}_{L^2(D)}$. The hypothesis
in~\eqref{eq:projected-chord-condition} therefore implies
$E_{\mathrm{in}}(f_\eta)>E_{\mathrm{out}}(f_\eta)$. Since $\eta\mapsto f_\eta$ is
linear and injective, Lemma~\ref{lem:energy-test} completes the proof.
\end{proof}

\subsection{The Dirichlet comparison and eigenvalue counting}

To prove Theorem~\ref{thm:counting}, recall that
\[
 0<\mu_1(D)\leq\mu_2(D)\leq\cdots
\]
are the Dirichlet eigenvalues on $D$, repeated according to multiplicity.  The
parameters $h_D$ and $g_D$ are those defined in~\eqref{geomsetup2}.

\begin{lemma}
\label{ljnkajdnakjnd}
Let $\Omega$ be a bounded Lipschitz domain satisfying the chord condition along
$\omega$. Let $D\subset \Omega(\omega)$ be a bounded Lipschitz domain with
$\overline D\subset \Omega(\omega)$, and assume that $\mathsf{a}$ and
$\mathsf{b}$ are Lipschitz near $\overline D$. If, for some $k\in\mathbb N$,
\begin{equation}
 \frac{4}{h_D}>g_D\mu_k(D)+2\sqrt{\mu_k(D)},
 \label{eq:quantitative-patch}
\end{equation}
then
\begin{equation}
 N_-(K^*)\geq k.
 \label{eq:counting-conclusion}
\end{equation}
In particular, we have the lower bound
\begin{equation}
  N_-(K^*)
  \geq
  N_D^{\mathrm{Dir}}\left(\frac{16}{(h_D+\sqrt{h_D}\sqrt{h_D+4g_D})^2}\right).
 \label{eq:dirbound}
\end{equation}
\end{lemma}

\begin{proof}
Let $\mathcal E_k$ be the span of the first $k$ Dirichlet eigenfunctions on
$D$.  For $\eta\in\mathcal E_k$, the variational characterization of
$\mu_k(D)$ gives
\begin{equation}
\mathfrak{E}_D(\eta)\leq\mu_k(D)\norm{\eta}_{L^2(D)}^2.
 \label{eq:dirichlet-comparison}
\end{equation}
Furthermore,
\begin{equation}
 4\int_D\frac{\eta^2}{\ell}\dd y
 \geq\frac{4}{h_D}\norm{\eta}_{L^2(D)}^2,
 \qquad
 \int_D(W_D-\ell)\abs{\nabla\eta}^2\dd y
 \leq g_D\mathfrak{E}_D(\eta).
 \label{eq:counting-upper}
\end{equation}
Combining~\eqref{eq:dirichlet-comparison}--\eqref{eq:counting-upper} with the
strict inequality~\eqref{eq:quantitative-patch} shows
that~\eqref{eq:projected-chord-condition} holds for every nonzero
$\eta\in\mathcal E_k$. The associated boundary data form a $k$-dimensional space,
and Lemma~\ref{lem:energy-test} gives~\eqref{eq:counting-conclusion}.

To obtain~\eqref{eq:dirbound}, note that the inequality
$4/h_D>g_D\Lambda+2\sqrt{\Lambda}$ is satisfied if and only if
\[
  \Lambda<\frac{16}{(h_D+\sqrt{h_D}\sqrt{h_D+4g_D})^2}.
\]
If $k\geq1$ is the value of the counting function on the right-hand side
of~\eqref{eq:dirbound}, then $\mu_k(D)$ lies below this threshold, so
that~\eqref{eq:quantitative-patch} holds for this $k$
and~\eqref{eq:counting-conclusion} gives~\eqref{eq:dirbound}; for $k=0$ there is
nothing to prove.
\end{proof}

We write
\[
 \mathfrak{W}(\omega,\Omega)
 \coloneqq
 \sup_{x\in\Omega}x\cdot\omega-\inf_{x\in\Omega}x\cdot\omega
\]
for the width of $\Omega$ in the direction $\omega$. We have $h_D\leq
\mathfrak{W}(\omega,\Omega)$ and $g_D\leq \mathfrak{W}(\omega,\Omega)$. Recall
the notation $N_D^{\mathrm{Dir}}(\Lambda)$ from~\eqref{defknefknd} for the
Dirichlet counting function of the projected patch $D$.

\begin{corollary}
\label{cor:flattening-count}
Let $\Omega$ be a bounded Lipschitz domain satisfying the chord condition
along~$\omega$. For every bounded Lipschitz domain $D\subset \Omega(\omega)$ with
$\overline D\subset \Omega(\omega)$ and $\mathsf{a},\mathsf{b}$ Lipschitz
near~$\overline D$,
\begin{equation}
 N_-(K^*)\geq
 N_D^{\mathrm{Dir}} \biggl(\frac{(\sqrt5-1)^2}{\mathfrak{W}(\omega,\Omega)^2}\biggr).
 \label{eq:flattening-count}
\end{equation}
Moreover, let $\Omega_\varepsilon$ be a family of bounded smooth domains
satisfying the chord condition along~$\omega$, and suppose that
$\mathfrak{W}(\omega,\Omega_\varepsilon)\to0$ as $\varepsilon\to0$. Let $D$ be a fixed
bounded Lipschitz domain with $\overline D\subset\Omega_\varepsilon(\omega)$ for
every $\varepsilon$, and assume that the endpoint functions are Lipschitz near
$\overline D$ for each $\varepsilon$. Then the family of NP operators
$K_\varepsilon^*$ on $\partial\Omega_\varepsilon$ satisfies
\begin{equation}
 N_-(K_\varepsilon^*)\geq
 \frac{\measure{D}(\sqrt5-1)^2}{4\pi \mathfrak{W}(\omega,\Omega_\varepsilon)^2}
 +o\!\left(\mathfrak{W}(\omega,\Omega_\varepsilon)^{-2}\right).
 \label{eq:flattening-weyl}
\end{equation}
\end{corollary}

\begin{proof}
Since $h_D\leq \mathfrak{W}(\omega,\Omega)$ and $g_D\leq
\mathfrak{W}(\omega,\Omega)$, the bound~\eqref{eq:flattening-count} follows
from~\eqref{eq:dirbound} in Lemma~\ref{ljnkajdnakjnd} and the fact that
$16/(1+\sqrt{5})^2=(\sqrt{5}-1)^2$.

For the family $\Omega_\varepsilon$, the Dirichlet Weyl law $
N_D^{\mathrm{Dir}}(\Lambda)=(\measure{D}/4\pi)\Lambda+o(\Lambda)$, as
$\Lambda\to\infty$, applied with $\Lambda =
(\sqrt5-1)^2/\mathfrak{W}(\omega,\Omega_\varepsilon)^2$,
gives~\eqref{eq:flattening-weyl}.
\end{proof}

For oblate spheroids~\eqref{eq:spheroidequation} with semi-axes $(a,a,b)$ this
estimate can be made explicit.

\begin{corollary}
\label{cor:oblate-flattening-count}
If $K_R^*$ denotes the NP operator on an oblate spheroid
with aspect ratio $R=a/b$, then
\begin{equation}
 \liminf_{R\to\infty}\frac{N_-(K_R^*)}{R^2}
 \geq\frac{(\sqrt5-1)^2}{16}.
 \label{eq:oblate-flattening-count}
\end{equation}
\end{corollary}

\begin{proof}
Let $R=a/b>1$ and
choose, for a fixed $0<\vartheta<1$, the projected disk
\[
 D_\vartheta=D(0,\vartheta a)=\set{y\in\R^2:\abs y<\vartheta a}
\]
in the $x_1x_2$-plane.  The width in the $x_3$-direction is $2b$, and
$\measure{D_\vartheta}=\pi\vartheta^2a^2$.
Corollary~\ref{cor:flattening-count} gives
\[
 N_-(K_R^*)\geq
 N_{D_\vartheta}^{\mathrm{Dir}}\biggl(\frac{(\sqrt5-1)^2}{4b^2}\biggr).
\]
The scaling law for the Dirichlet spectrum of a disk and the planar Weyl law
therefore yield
\[
 N_{D_\vartheta}^{\mathrm{Dir}}\biggl(\frac{(\sqrt5-1)^2}{4b^2}\biggr)
 =\frac{\vartheta^2(\sqrt5-1)^2}{16}R^2+o(R^2).
\]
Taking the limit inferior and then letting $\vartheta\uparrow1$
proves~\eqref{eq:oblate-flattening-count}.
\end{proof}

We record two consequences of Corollary~\ref{cor:flattening-count} with
elementary hypotheses.

\begin{example}
\label{ex:disk}
 If $D=D(0,\rho)$ is chosen to be a disk of radius $\rho$
contained in $\Omega(\omega)$, then
\[
 \mu_1(D)=\frac{j_{0,1}^2}{\rho^2},
\]
where $j_{0,1}$ is the first positive zero of the Bessel function $J_0$. Hence if
$\Omega$ satisfies the chord condition along $\omega$ and
$A(\overline{D(0,\rho)})\subset \Omega(\omega)$ for some affine isometry $A$,
then
\[
 \frac{\mathfrak{W}(\omega,\Omega)}{\rho}<\frac{\sqrt{5}-1}{j_{0,1}}
 \implies
 N_-(K^*)\geq1.
\]
In other words, we can ensure that $N_-(K^*)\geq 1$ as soon as $\Omega$ satisfies
the chord condition along some $\omega$ such that there is a closed disk of
radius $\rho>j_{0,1}\mathfrak{W}(\omega,\Omega)/(\sqrt{5}-1)$ in the interior of
$\Omega(\omega)$. The constant $(\sqrt{5}-1)/j_{0,1}$ is approximately $0.514$.
Replacing the first eigenfunction by $\eta_\rho(y)\coloneqq 1-\abs y^2/\rho^2$
gives the elementary, slightly weaker condition
\[
 \frac{\mathfrak{W}(\omega,\Omega)}{\rho}<\frac{\sqrt{5}-1}{\sqrt6}.
\]
\end{example}

\begin{example}
For a rectangle $D=(-r_1,r_1)\times(-r_2,r_2)$, we have
\[
 \mu_1(D)=\frac{\pi^2}{4}\left(\frac1{r_1^2}+\frac1{r_2^2}\right).
\]
The width bound then reads as follows: if $\Omega$ satisfies the chord
condition along $\omega$ and $A([-r_1,r_1]\times[-r_2,r_2])\subset \Omega(\omega)$
for some affine isometry $A$, then
\[
 \mathfrak{W}(\omega,\Omega)<\frac{2(\sqrt5-1)r_1r_2}
  {\pi\sqrt{r_1^2+r_2^2}} \implies
 N_-(K^*)\geq1.
\]
This bound requires both tangential dimensions to be large compared with the
width; it then guarantees negative spectrum, and a negative
eigenvalue in the compact case. The bound distinguishes plates from rods: both
$r_1$ and $r_2$ must become large compared with $\mathfrak{W}(\omega,\Omega)$.
\end{example}

\begin{example}
\label{ex:oblate-chord}

For a surface of revolution about the $x_3$-axis, write the upper and lower parts
of the generating curve as $z_+(r)$ and $z_-(r)$, defined for $0\leq r\leq\rho$.
With our earlier notation, $\mathsf{b}(y)=z_+(\abs y)$ and
$\mathsf{a}(y)=z_-(\abs y)$.
Put
\[
 \ell(r)=z_+(r)-z_-(r),
 \qquad
 W_D=\max z_+-\min z_-.
\]
For a radial test function $\eta$ with $\eta(\rho)=0$,
condition~\eqref{eq:projected-chord-condition} becomes, after division by the
common angular factor,
\[
 4\int_0^\rho\frac{\eta(r)^2}{\ell(r)}r\dd r
 >
 \int_0^\rho\bigl(W_D-\ell(r)\bigr)\eta'(r)^2r\dd r
 +2\sqrt{\left(\int_0^\rho\eta(r)^2r\dd r\right)
 \left(\int_0^\rho\eta'(r)^2r\dd r\right)}.
\]

Consider the interior of the oblate spheroid
\[
 \Omega
 =
 \set[\bigg]{(x_1,x_2,x_3)\in \R^3:\frac{x_1^2+x_2^2}{a^2}+\frac{x_3^2}{b^2}<1},
 \qquad a>b.
\]
With $\omega$ in the $x_3$-direction,
\[
 \ell(r)=2b\sqrt{1-\frac{r^2}{a^2}},
 \qquad W_D=2b.
\]
Here we take $\eta=\eta_a$ from Example~\ref{ex:disk} on $D=\Omega(\omega)$, the
full disk of radius $a$, so the strict inclusion $\overline D\subset
\Omega(\omega)$ in Proposition~\ref{prop:projected-chord} is not satisfied. This
endpoint case is nevertheless admissible: on the spheroid, the resulting boundary
data are
\begin{equation}
 f_\eta(x_1,x_2,x_3)=\frac{x_3\abs{x_3}}{b^2},
 \label{eq:oblate-boundary-data}
\end{equation}
which are $C^{1,1}$ across the equator and vanish there. Alternatively, for
$0<\vartheta<1$ one may apply the proposition on the disk
$D_\vartheta=D(0,\vartheta a)$ with $\eta_{\vartheta a}$, extend $\eta_{\vartheta
a}$ by zero, and let $\vartheta\uparrow 1$. The zero-extended functions converge
to $\eta_a$ in $H^1(\Omega(\omega))$. The unweighted $L^2$ and gradient terms
then converge by dominated convergence. For the remaining term, the integrands
are dominated by $\eta^2/\ell$; near the rim, with $s=a-r$, this function is
$O(s^{3/2})$, so the singular weight $1/\ell$ is integrable. Consequently, all
three energy terms converge to the endpoint values below. Whenever the limiting
inequality is strict, it therefore also holds for $\vartheta<1$ sufficiently
close to $1$. Direct calculation gives
\[
 4\int_D\frac{\eta^2}{\ell}\dd y
 =\frac{4\pi a^2}{5b},
 \qquad
 \int_D(W_D-\ell)\abs{\nabla\eta}^2\dd y
 =\frac{28\pi b}{15},
 \qquad
 2\sqrt{\mathfrak{E}_D(\eta)}\norm{\eta}_{L^2(D)}
 =2\pi a\sqrt{\frac23}.
\]
It follows from the argument proving Proposition~\ref{prop:projected-chord} that
the condition
\begin{equation}
 \frac ab>\frac{7}{\sqrt6}\approx2.86
 \label{eq:oblate-aspect-ratio-condition}
\end{equation}
implies the existence of a negative eigenvalue of $K^*$ for the oblate
spheroid~\eqref{eq:spheroidequation} with semi-axes $(a,a,b)$.
\end{example}

\begin{example}
\label{ex:superellipsoid-cap}
Let $\Omega$ be the interior of the superellipsoid,
\[
 \Omega=\set{(x_1,x_2,x_3)\in\R^3:(x_1^2+x_2^2)^{k/2}+\abs {x_3}^k<1}.
\]
We show that $N_-(K^*)\geq1$ for $k>4$ using Lemma~\ref{lem:energy-test}
directly. 

Following the notation in the Introduction, 
we take $\omega=(0,0,1)$ and use the coordinates $y=(x_1,x_2)$ and $t=x_3$. 
Put $z_+(r)=(1-r^k)^{1/k}$ and, for $0<\rho<1$, take $D=D(0,\rho)$ and
$\eta=\eta_\rho$ from Example~\ref{ex:disk}, extended by zero. We define the
function $\tilde{f}_\eta\in H^{1/2}(\partial \Omega)$ by
$\tilde{f}_\eta(y,z_+(\abs y))=\eta(y)$ on the upper part of $\partial\Omega$ and
zero elsewhere on $\partial\Omega$. Unlike the functions $f_\eta$ in
Proposition~\ref{prop:projected-chord} and Example~\ref{ex:oblate-chord}, which
equal $\pm\eta$ on the upper and lower parts of $\partial\Omega$, the function
$\tilde f_\eta$ is supported on the upper part only.

As in the proof of Lemma~\ref{prop:flat-face}, let $u_H^-$ and $u_H^+$ be the
harmonic extensions of $\eta$ to $H_-=\set{t<1}$ and $H_+=\set{t>1}$, with
common energy $E_H(\eta)$. Extending $\eta$ constantly through the columns
$z_+(\abs y)<t<1$ and using $u_H^+$ above $t=1$ gives
\[
 E_{\mathrm{out}}(\tilde{f}_\eta)\leq E_H(\eta)+A,
 \qquad
 A=\int_D(1-z_+(\abs y))\abs{\nabla\eta}^2\dd y
 \leq\frac{8\pi\rho^k}{k(k+4)(1-\rho^k)}.
\]
The last estimate follows from
$1-(1-r^k)^{1/k}\leq r^k/[k(1-\rho^k)]$ for $r\leq\rho$.

For the interior comparison, use the quarter-space $Q=\set{y_1<1,\ t<1}$
containing $\Omega$, and reflect by $\mathsf{R}y=(2-y_1,y_2)$. The function
$u_H^-(y,t)-u_H^-(\mathsf{R}y,t)$ minimizes the energy in $Q$ with upper trace
$\eta$ and zero trace on the vertical face. The half-space Poisson kernel
$P_s(y)=s/[2\pi(\abs y^2+s^2)^{3/2}]$ gives its energy as $E_H(\eta)+B$, where
\[
 B=\frac1{2\pi}\int_D\int_D
 \frac{\eta(y)\eta(y')}{\abs{y-\mathsf{R}y'}^3}\dd y'\dd y
 \geq\frac{\pi\rho^4}{64(1+\rho)^3}.
\]
Indeed, the interaction energy of $u_H^-$ and its reflection is $-B$; the
inequality uses $\int_D\eta\dd y=\pi\rho^2/2$ and
$\abs{y-\mathsf{R}y'}\leq2(1+\rho)$. Extending $u_{\mathrm{in}}$ to $Q$ by the
same constant columns and by zero elsewhere gives
$E_{\mathrm{in}}(\tilde{f}_\eta)+A\geq E_H(\eta)+B$.
Consequently,
\[
 E_{\mathrm{in}}(\tilde{f}_\eta)-E_{\mathrm{out}}(\tilde{f}_\eta)
 \geq B-2A
 \geq\frac{\pi\rho^4}{64(1+\rho)^3}
 -\frac{16\pi\rho^k}{k(k+4)(1-\rho^k)},
\]
which is positive for $\rho$ sufficiently small since $k>4$.
Lemma~\ref{lem:energy-test} therefore gives a negative eigenvalue.
\end{example}

\section{A balanced-Hessian condition excluding negative eigenvalues}
\label{sec:positivity}

The results of Section~\ref{sec:chords} can be used to detect negative spectrum
of the NP operator. We next derive a condition that excludes negative eigenvalues.
The condition is expressed through a defining function
$F$ and follows from a Rellich identity for the single-layer potential. The
required hypothesis $(\Delta F)I\geq 2D^2 F$ is a condition on the defining
function $F$ in the ambient space $\R^3$, not solely a condition on boundary
curvatures. Since we are in dimension three, the condition $(\Delta F)I\geq 2D^2
F$ implies convexity of $F$: if $\nu_1\leq\nu_2\leq\nu_3$ are the Hessian
eigenvalues, then $\nu_1+\nu_2\geq\nu_3$ implies $\nu_1\geq0$. We now prove
Theorem~\ref{thm:balanced-hessian} on the absence of negative eigenvalues.

\begin{proof}[Proof of Theorem~\ref{thm:balanced-hessian}]
Suppose that $K^*\varphi=\lambda\varphi$ for some nonzero
$\varphi\in\mathfrak{H}$ and $\lambda\in \R\setminus\set{0}$; the case
$\lambda=0$ is immediate. We show that $\lambda\geq0$. The integral
kernel of $K^*$ is real, so by replacing $\varphi$ by its
nonzero real or imaginary part if necessary, we may assume that $\varphi$ is
real-valued. Since $0$ is a regular value of $F$, the boundary is $C^2$. The
compact realizations of $K^*$ on $L^2(\partial\Omega)$ and
$\mathfrak{H}$ have the same nonzero eigenspaces: this follows by
restricting their finite-rank Riesz projections to the dense subspace
$L^2(\partial\Omega)$. Thus $\varphi\in L^2(\partial\Omega)$. Put $u=\mathcal S
\varphi$. Standard single-layer estimates give $u|_{\partial\Omega}\in
H^1(\partial\Omega)$ and one-sided normal derivatives in $L^2(\partial\Omega)$;
see~\cite[Section~1.7]{ChandlerWildeSpence2024}. Moreover, $u=O(\abs{x}^{-1})$ and
$\nabla u=O(\abs{x}^{-2})$ at infinity.

Write $q_{\mathrm{in}}$ and $q_{\mathrm{out}}$ for the normal derivatives of $u$
on the two sides of the boundary, both in the direction of the outward normal $n$
of $\Omega$. The jump relations~\eqref{eq:jump-in-out} give
\begin{equation}
\label{eq:jhumprel}
 q_{\mathrm{out}}-q_{\mathrm{in}}=\varphi,
 \qquad
 q_{\mathrm{out}}+q_{\mathrm{in}}=K^*\varphi.
\end{equation}

We call
\[
 \mathcal R(v,F)
 \coloneqq\abs{\nabla v}^2\nabla F-2(\nabla F\cdot\nabla v)\nabla v
\]
the Rellich vector field. For a harmonic function $v$, direct differentiation
gives
\begin{equation}
\label{reladlad}
 \operatorname{div}\mathcal R(v,F)
 =(\Delta F)\abs{\nabla v}^2-2D^2F[\nabla v,\nabla v].
\end{equation}
Since $F$ is constant on $\partial\Omega$ and $\Omega=\set{x: F(x)<0}$, we have
$\nabla F=\abs{\nabla F}\,n$ on $\partial\Omega$. When~\eqref{reladlad} is
applied inside and outside $\Omega$, the tangential boundary terms cancel because
the boundary trace of $u$ is continuous. Using~\eqref{eq:jhumprel} for the normal
terms gives
\begin{align}
  \int_{\partial\Omega}\abs{\nabla F}(K^*\varphi)\varphi\dd S
 ={}&\int_\Omega\left((\Delta F)\abs{\nabla u}^2
 -2D^2F[\nabla u,\nabla u]\right)\dd x
 \notag\\
 &+\int_{\R^3\setminus\overline\Omega}
 \left((\Delta F)\abs{\nabla u}^2
 -2D^2F[\nabla u,\nabla u]\right)\dd x.
 \label{eq:np-rellich-identity}
\end{align}
The integrated Rellich identity is valid with these $L^2$ boundary traces
by~\cite[Lemmas~4.6--4.7]{ChandlerWildeSpence2024}, applied with the multiplier
$\nabla F$ on bounded truncations. The bulk integrals converge because $D^2F$ is
bounded and $u$ has finite Dirichlet energy. Boundedness of $D^2F$ also gives
$\nabla F=O(1+\abs{x})$, so the boundary term on a sphere of radius $R$ is
$O(R^{-1})$ and tends to zero as $R\to \infty$. Since $(\Delta F)I-2D^2F\geq0$,
the right-hand side of~\eqref{eq:np-rellich-identity} is nonnegative. The
left-hand side equals
\[
 \lambda\int_{\partial\Omega}\abs{\nabla F}\abs{\varphi}^2\dd S,
\]
and the integral is positive for $\varphi\neq0$, so $\lambda\geq 0$.
\end{proof}

\subsection{The squared gauge and ellipsoids}

There is a natural defining function associated with a convex domain. We can assume
that the origin lies inside $\Omega$ and let
\[
 p_\Omega(x)\coloneqq\inf\set{t>0:x\in t\Omega}
\]
be its Minkowski gauge.  Then $p_\Omega=1$ on $\partial\Omega$ and
$p_\Omega(tx)=tp_\Omega(x)$ for $t>0$.  Define
\[
 F_\Omega(x)=\frac12p_\Omega(x)^2.
\]
We have $\Omega=\set{F_\Omega<1/2}$; subtracting $1/2$ changes neither the gradient
nor the Hessian. The function $F_\Omega$ is homogeneous of degree two, so
$D^2F_\Omega$ is homogeneous of degree zero. It therefore suffices to check the
Hessian on a single transversal surface.

\begin{corollary}
\label{cor:squared-gauge}
Let $\Omega$ be a smooth convex domain with $0\in\Omega$, and
let $F_\Omega$ be the squared gauge defined above. If
$(\Delta F_\Omega)I-2D^2F_\Omega\geq0$ on $\partial\Omega$, then $K^*$ has no
negative eigenvalues.
\end{corollary}

\begin{proof}
The function $F_\Omega$ is smooth away from the origin. By homogeneity, its
Hessian is bounded and the assumed inequality holds throughout
$\R^3\setminus\set{0}$. Euler's identity gives
$x\cdot\nabla F_\Omega=2F_\Omega=1$ on $\partial\Omega$, so $\nabla F_\Omega$ does
not vanish there.
The argument in the proof of Theorem~\ref{thm:balanced-hessian} applies after
removing a small ball around the origin.  To see that the new boundary term
disappears, consider the Rellich vector field $\mathcal R(u,F_\Omega)$
from~\eqref{reladlad}, with $u=\mathcal S \varphi$ as in that proof.  Since
$F_\Omega$ is homogeneous of degree two, $\nabla F_\Omega(x)=O(\abs{x})$ as
$x\to0$, whereas $u$, being harmonic in $\Omega$, has bounded gradient near the
origin.  Hence
$\mathcal R(u,F_\Omega)\cdot n=O(r)$ on $\partial B_r$, and
\[
 \int_{\partial B_r}\mathcal R(u,F_\Omega)\cdot n\,\dd S=O(r^3)\to 0
 \qquad (r\to 0).
\]
\end{proof}

\begin{example}\label{cor:ellipsoid-balanced}
Consider the ellipsoid with semi-axes $a_1\geq a_2\geq a_3>0$,
\[
\set[\bigg]{x=(x_1,x_2,x_3)\in\R^3:
 \frac{x_1^2}{a_1^2}+\frac{x_2^2}{a_2^2}+\frac{x_3^2}{a_3^2}=1}.
\]
For the quadratic defining function
\[
 F(x)=\frac12\left(\frac{x_1^2}{a_1^2}+
 \frac{x_2^2}{a_2^2}+\frac{x_3^2}{a_3^2}-1\right),
\]
the Hessian eigenvalues are $a_1^{-2},a_2^{-2},a_3^{-2}$. We obtain the following
condition: if
\begin{equation}
 \frac1{a_3^2}\leq\frac1{a_1^2}+\frac1{a_2^2},
 \label{eq:ellipsoid-balanced-condition}
\end{equation}
then the NP operator $K^*$ has no negative eigenvalues.

In other words, the reciprocal squares of the semi-axes satisfy the triangle
inequality.  Equivalently,
\[
 a_3\geq\frac{a_1a_2}{\sqrt{a_1^2+a_2^2}}.
\]

For a prolate spheroid with semi-axes $(a,a,b)$, $b>a$,
condition~\eqref{eq:ellipsoid-balanced-condition} is automatic. Thus every
prolate spheroid satisfies the balanced-Hessian condition, independently of its
aspect ratio. For an oblate spheroid with semi-axes $(a,a,b)$, the condition
becomes
\[
 \frac ab\leq\sqrt2.
\]
This gives a simple region of guaranteed nonnegativity, in contrast with
condition~\eqref{eq:oblate-aspect-ratio-condition} in
Example~\ref{ex:oblate-chord}.
\end{example}

\subsection{Stability under perturbations}

The balanced-Hessian condition $(\Delta F)I\geq 2D^2 F$ is stable under a
quantitative $C^2$ perturbation of a quadratic defining function. Let $A$ be
positive definite with eigenvalues $0<\alpha_1\leq\alpha_2\leq\alpha_3$ and set
\[
 \delta_A\coloneqq\alpha_1+\alpha_2-\alpha_3.
\]

\begin{proposition}
\label{prop:balanced-perturbation}
Let $F(x)=\tfrac12x^TAx+g(x)$ with $g\in C^2(\R^3)$. Assume that $\delta_A>0$.  If
\[
 \sup_{x\in\R^3}\norm{D^2g(x)}_{\mathrm{op}}<\frac{\delta_A}{3},
\]
then the NP operator on the boundary of every bounded regular level surface
$\set{x\in \R^3 : F(x)=c}$ has no negative eigenvalues.
\end{proposition}

\begin{proof}
Put $\varepsilon=\sup_x\norm{D^2g(x)}_{\mathrm{op}}$.  If
$\nu_1\leq\nu_2\leq\nu_3$ are the eigenvalues of $D^2F$, Weyl's inequalities
give
\[
 \nu_3\leq\alpha_3+\varepsilon,
 \qquad
 \nu_1+\nu_2\geq\alpha_1+\alpha_2-2\varepsilon.
\]
Since $\delta_A\leq\alpha_1$, the same hypothesis gives
$\varepsilon<\alpha_1$, and therefore
$\nu_1\geq\alpha_1-\varepsilon>0$.  Thus $D^2F$ is positive definite, so
$F$ is convex.  Moreover, $3\varepsilon<\delta_A$ gives
$\nu_3<\nu_1+\nu_2$.  The bounded-Hessian assumption also implies the
required linear growth of $\nabla F$.  Hence
Theorem~\ref{thm:balanced-hessian} applies.
\end{proof}

For the unit ball, $A=I$ and $\delta_A=1$.  Thus
\[
 \sup_x\norm{D^2g(x)}_{\mathrm{op}}<\frac13
\]
guarantees that the NP operator is nonnegative for every bounded regular level
surface of $F(x)=\tfrac12\abs x^2+g(x)$. As an explicit strongly convex,
nonellipsoidal example, set
\begin{equation}
 F_\varepsilon(x)=\frac12\abs x^2+\varepsilon\log\cosh x_3,
 \qquad 0 < \varepsilon\leq1.
 \label{eq:logcosh-body}
\end{equation}
Then
\[
 D^2F_\varepsilon=\operatorname{diag}\bigl(1,1,
 1+\varepsilon\sech^2x_3\bigr),
\]
so the largest eigenvalue is at most $2$, while the other two sum to $2$.
Hence Theorem~\ref{thm:balanced-hessian} applies directly (the perturbation
$g(x)=\varepsilon\log\cosh x_3$ satisfies
$\sup_x\norm{D^2g(x)}_{\mathrm{op}}=\varepsilon$, which may exceed $1/3$), and
the NP operator on every bounded regular level surface of~\eqref{eq:logcosh-body}
has no negative eigenvalues.

\section{Numerical motivation for theoretical results}
\label{sec:numerics}

Numerical computations on surfaces of revolution motivated our theoretical
criteria for the presence of negative eigenvalues. We now present (mostly)
graphical results from some of these computations. The eigenvalues in the figures
below are computed using the high-order axisymmetric Fourier--Nyström
framework of~\cite[Section~9.2]{HKR} (based on~\cite{YoungHaoMartinsson}) and
matrix eigensolvers in \textsc{Matlab}. The key idea in the Fourier--Nyström
framework is the decomposition into Fourier modes in the azimuthal variable and a
one-dimensional discretization along the generating curve. See
Figure~\ref{fig:curves} for an overview of the surfaces of revolution we
consider. The spheroids serve as smooth calibration examples, the cheese and the
Swedish coconut macaroon illustrate Theorem~\ref{thm:flat}, and the sausage and
superellipsoid probe other geometric transitions.

All displayed eigenvalues are computed approximations of eigenvalues of $K^*$ in
the normalization of~\eqref{eq:intro-np}, with spectrum in $[-1,1]$ and
distinguished eigenvalue $1$. In the full spectrum, each eigenvalue for
azimuthal index $m>0$ has azimuthal multiplicity two, whereas an eigenvalue for
$m=0$ has azimuthal multiplicity one.

The one-dimensional Nyström discretization produces a finite matrix whose
dimension is adjusted to the relevant geometric parameter range. The matrices we
use have size at most $3072\times3072$. We compute the relevant matrix
eigenvalues and plot selected ones as curves. We visualize around 120 eigenvalue
branches for each example, usually associated with the largest eigenvalues (in
modulus). Selection by modulus need not follow a single eigenvalue branch through
a branch crossing. The largest eigenvalues are well resolved, whereas the
smallest (in modulus) computed eigenvalues are typically sensitive to the
discretization and need not correspond to a particular exact eigenvalue.

\subsection{Spheroids}

\begin{figure}[htbp]
    \centering
    \begin{subfigure}{0.48\textwidth}
        \centering
        \includegraphics[width=\linewidth]{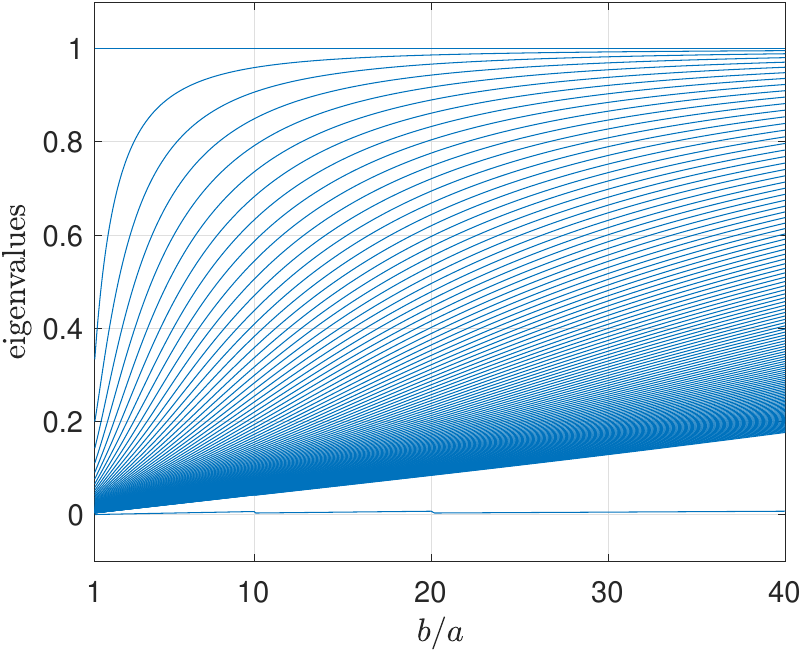}
        \caption{$m=0$}
        \label{fig:prolate0}
    \end{subfigure}
    \hfill
    \begin{subfigure}{0.48\textwidth}
        \centering
        \includegraphics[width=\linewidth]{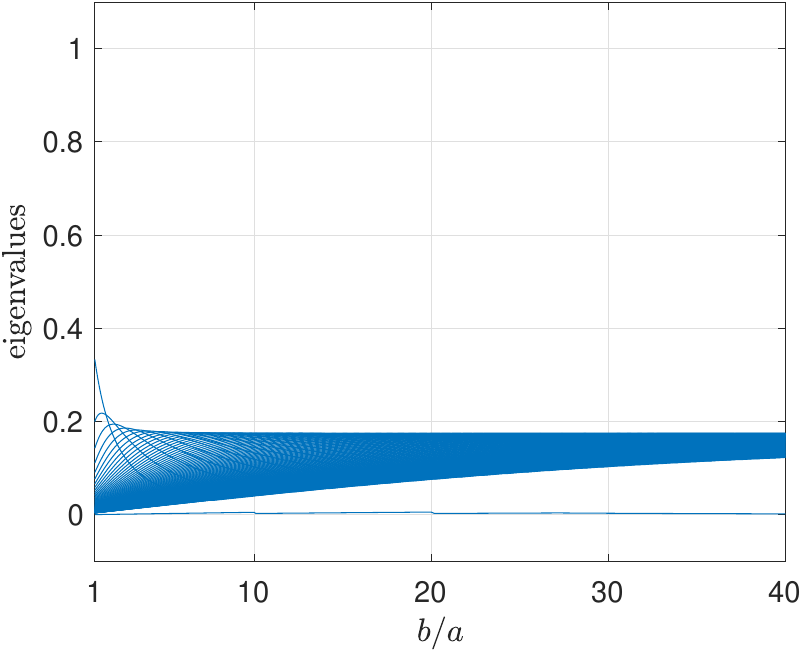}
        \caption{$m=1$}
        \label{fig:prolate1}
    \end{subfigure}
    \caption{Computed eigenvalues for the prolate spheroid as a function of
    $b/a\geq1$. The 120 largest computed eigenvalues and the smallest computed
    eigenvalue are shown for each azimuthal index.}
    \label{fig:prolate}
\end{figure}

The prolate spheroid is given by
\[
 \frac{x_1^2+x_2^2}{a^2}+\frac{x_3^2}{b^2}=1,
 \qquad b>a.
\]
The computations for azimuthal indices $m=0$ and $m=1$ detect no negative
eigenvalues as $b/a$ increases; see Figure~\ref{fig:prolate}. This agrees with
Example~\ref{cor:ellipsoid-balanced}, where we conclude nonnegativity for every
prolate spheroid. A large aspect ratio alone does not force a negative
eigenvalue.

\begin{figure}[htbp]
    \centering
    \begin{subfigure}{0.48\textwidth}
        \centering
        \includegraphics[width=\linewidth]{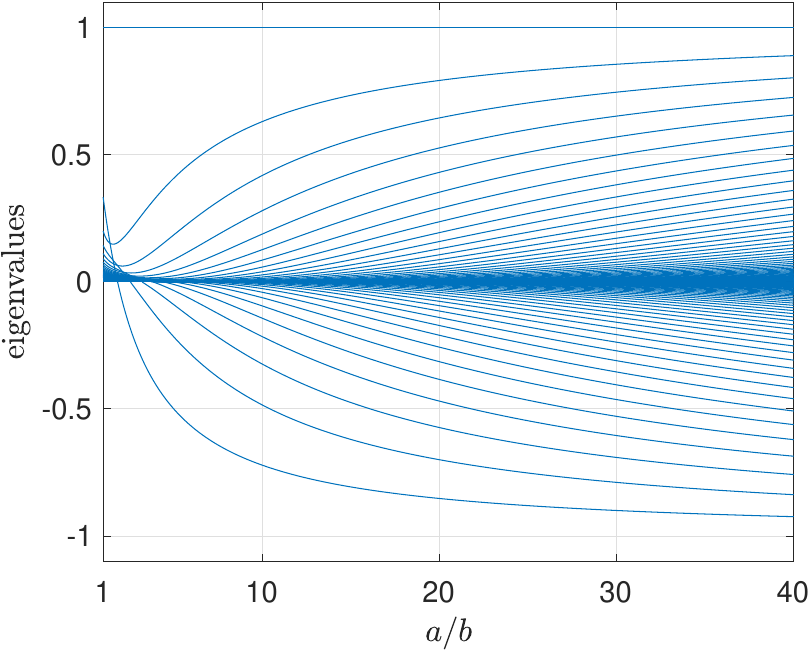}
        \caption{$m=0$}
        \label{fig:oblate0}
    \end{subfigure}
    \hfill
    \begin{subfigure}{0.48\textwidth}
        \centering
        \includegraphics[width=\linewidth]{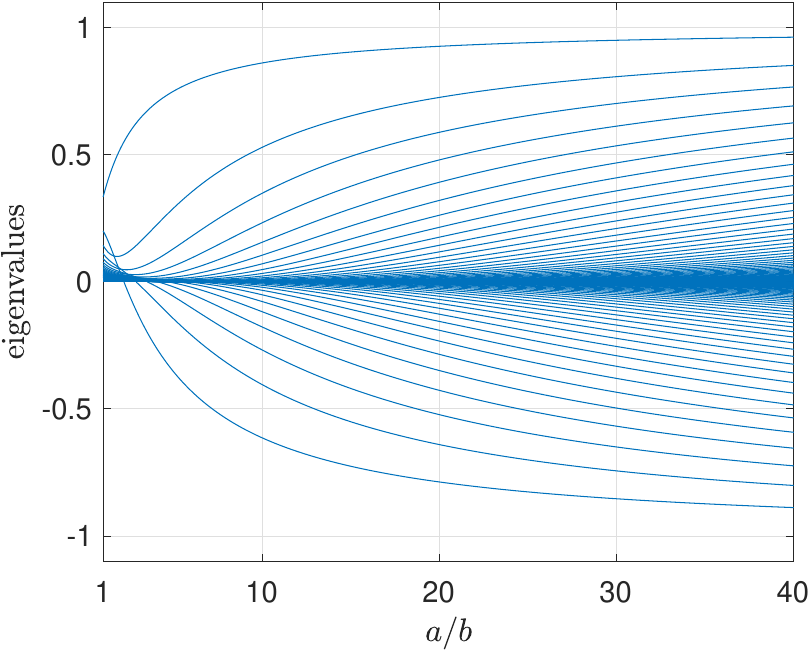}
        \caption{$m=1$}
        \label{fig:oblate1}
    \end{subfigure}
    \caption{Computed eigenvalues for the oblate spheroid as a function of
    $a/b\geq1$. The 120 computed eigenvalues of largest modulus are shown for
    each azimuthal index.}
    \label{fig:oblate}
\end{figure}

\begin{figure}[htb]
  \centering
  \includegraphics[height=5.75cm]{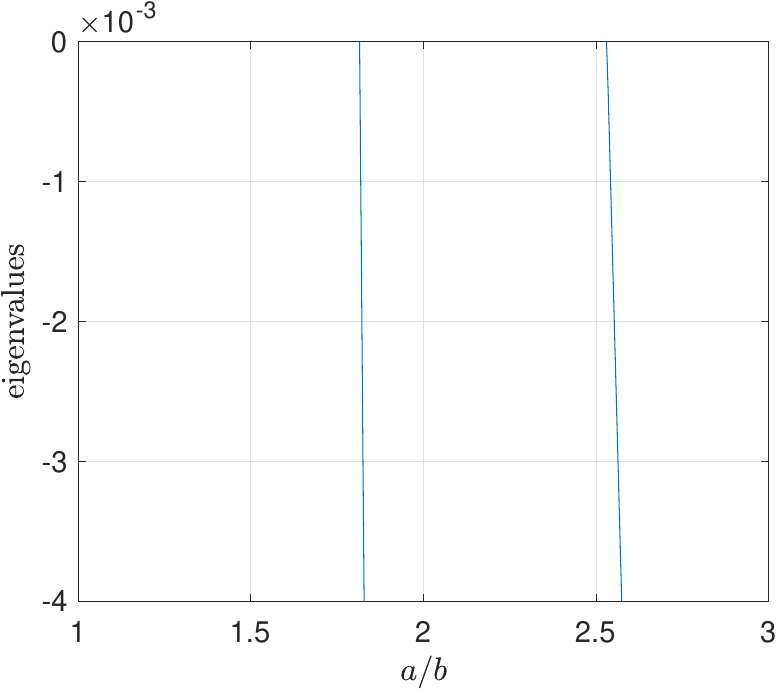}
  \caption{A close-up of the curves in
  Figure~\ref{fig:oblate}(\subref{fig:oblate0}) showing the first two zero
  crossings for the oblate spheroid. The value $R_*\approx 1.81527095947903$ of
  the first zero crossing computed numerically is accurate almost to machine
  precision; see~\eqref{eq:Rstar}.}
  \label{fig:oblate0zoom}
\end{figure}

For the oblate spheroid
\[
 \frac{x_1^2+x_2^2}{a^2}+\frac{x_3^2}{b^2}=1,
 \qquad a>b,
\]
some of the computed branches become negative as $a/b$ increases; see
Figure~\ref{fig:oblate}. In Figure~\ref{fig:oblate0zoom} we zoom in to see that
the first zero crossing found numerically is consistent with $R_*$
in~\eqref{eq:Rstar}.

Figure~\ref{fig:oblate} shows that the more oblate the spheroid becomes, the more
negative eigenvalues appear. On the other hand, since oblate spheroids have
everywhere positive principal curvatures, for a fixed aspect ratio $a/b$ there
exist at most finitely many negative eigenvalues~\cite[Theorem~3]{MR3D}. This
distinction between finitely many negative eigenvalues at fixed $a/b$ and an
increasing number as $a/b\to\infty$ is also emphasized in~\cite{AKLM22}.

\subsection{Comparison with the exact oblate formulas}
\label{sec:oblate-exact}

\begin{figure}[htbp]
    \centering
        \includegraphics[height=6cm]{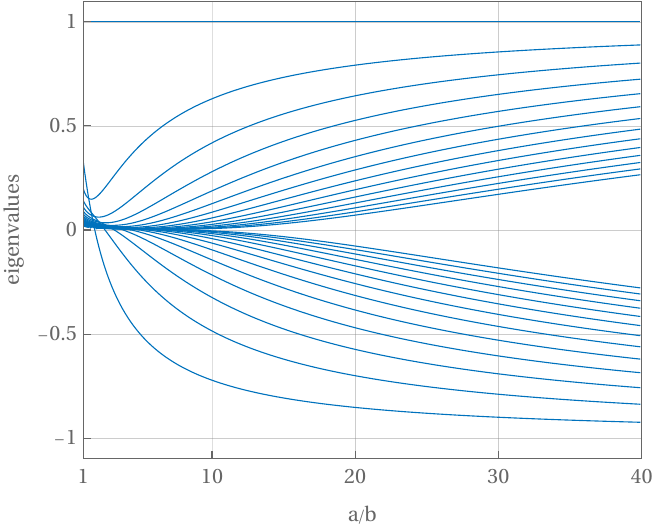}
    \caption{Eigenvalues $\lambda_{0,n}$ obtained from the exact
    formula~\eqref{eq:oblate-ahner-formula} for the oblate spheroid with $m=0$
    and $0\leq n\leq 26$.}
    \label{fig:comparisons}
\end{figure}

For completeness, we compare the high-order discretization with the exact
formulas for oblate spheroids derived by Ahner~\cite{Ahner}, building on the
earlier spheroidal analysis of Ahner and Arenstorf~\cite{AhnerArenstorf86}.
Let
\[
 \xi=\frac{1}{\sqrt{(a/b)^2-1}},
\]
and let $P_n^m$ and $Q_n^m$ denote the associated Legendre functions on the
principal branches of~\cite[Section~14.21]{DLMF}, continued from real arguments
greater than $1$ through the upper half-plane to $\ii\xi$. Note that this branch
of $Q_n^m$ decays at infinity. In the normalization used here, the eigenvalues of
$K^*$ are (here we reparametrize the eigenvalues from the introduction and use
$\xi=(R^2-1)^{-1/2}$ as the parameter)
\begin{equation}\label{eq:oblate-ahner-formula}
 \lambda_{m,n}(\xi)
 =1-2\frac{(n-m)!}{(n+m)!}(\xi^2+1)
 \bigl[\ii^{n+2m+1}Q_n^m(\ii\xi)\bigr]
 \biggl[\ii^{-n}\frac{\dd}{\dd\xi}P_n^m(\ii\xi)\biggr],
\end{equation}
where $0\leq m\leq n$. Ahner uses the opposite sign convention with
the same normalization, so the expression in~\eqref{eq:oblate-ahner-formula} is
the negative of his eigenvalue formula. In particular, it gives
$\lambda_{0,0}=1$, as required by~\eqref{eq:intro-np}.
Figure~\ref{fig:comparisons} shows eigenvalues
calculated from the exact formula~\eqref{eq:oblate-ahner-formula}.
The values are essentially the same as those obtained by purely numerical computations
in Figure~\ref{fig:oblate}(\subref{fig:oblate0}).

\subsection{The cheese, sausage, and Swedish coconut macaroon}

The remaining convex examples investigate the role of flat faces and tangential
junctions. These surfaces are $C^{1,1}$: regular enough for
Theorems~\ref{thm:flat} and~\ref{thm:counting}, but not $C^2$ as
Theorem~\ref{thm:balanced-hessian} requires, since the curvature jumps. The
cheese and the Swedish coconut macaroon have flat faces, so
Theorem~\ref{thm:flat} gives infinitely many negative eigenvalues. The sausage
has no flat face.

\begin{figure}[htbp]
    \centering
    \begin{subfigure}{0.48\textwidth}
        \centering
        \includegraphics[width=\linewidth]{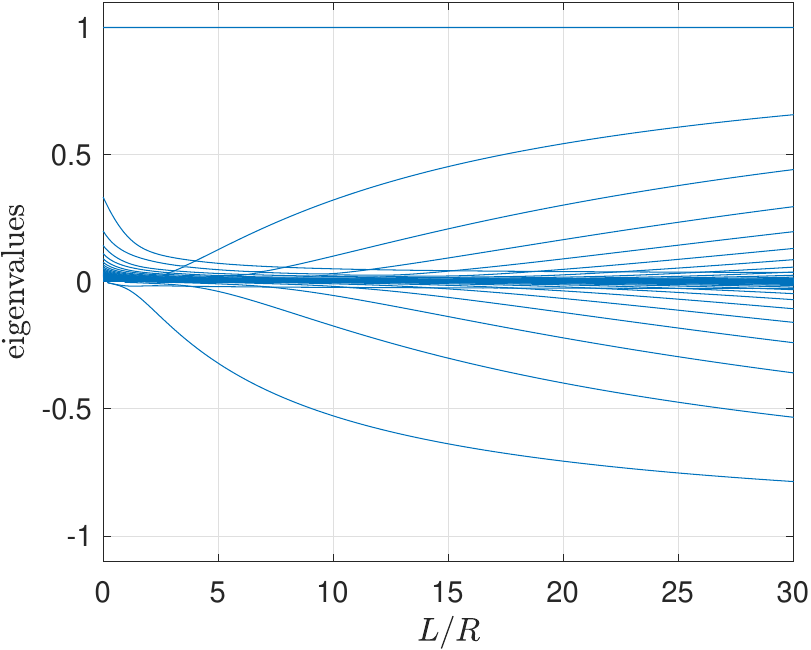}
        \caption{$m=0$}
        \label{fig:ost0}
    \end{subfigure}
    \hfill
    \begin{subfigure}{0.48\textwidth}
        \centering
        \includegraphics[width=\linewidth]{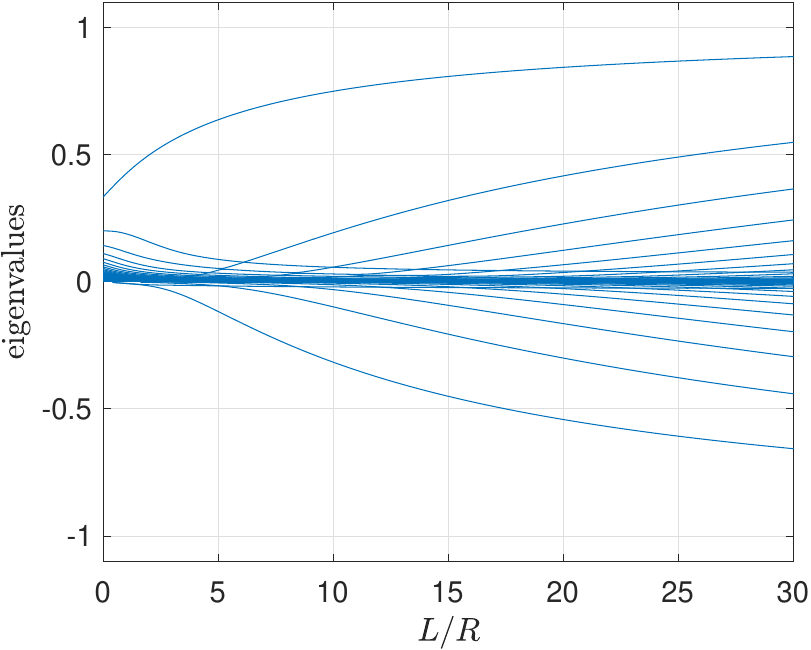}
        \caption{$m=1$}
        \label{fig:ost1}
    \end{subfigure}
    \caption{The cheese: computed eigenvalues as a function of $L/R$, where $R$
    is the radius of the semicircles and $L$ the length of the straight
    segments.}
    \label{fig:cheese}
\end{figure}

The cheese profile is obtained by separating two semicircles of radius $R$ and
joining them by straight segments of length $L$. The resulting surface has a flat
face for every positive separation, so infinitely many negative eigenvalues
appear as soon as $L>0$; see Figure~\ref{fig:cheese}.

\begin{figure}[htbp]
    \centering
    \begin{subfigure}{0.48\textwidth}
        \centering
        \includegraphics[width=\linewidth]{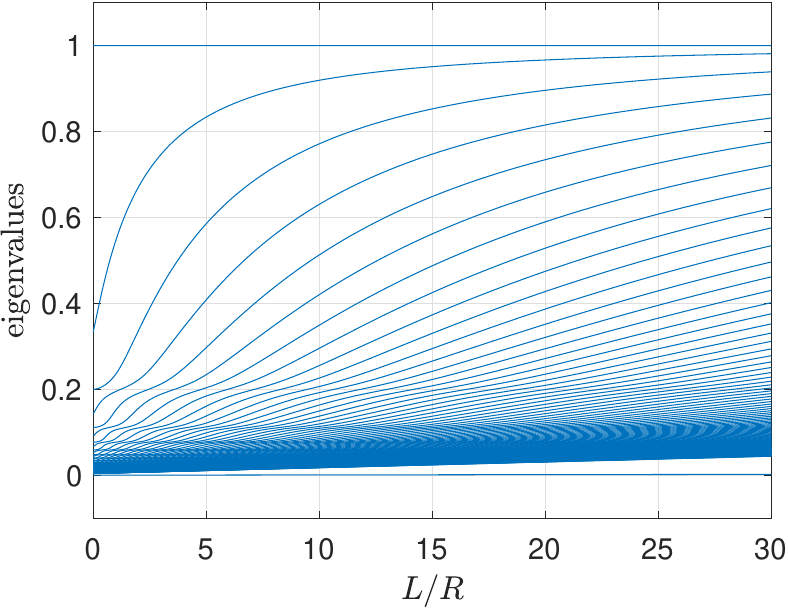}
        \caption{$m=0$}
        \label{fig:korv0}
    \end{subfigure}
    \hfill
    \begin{subfigure}{0.48\textwidth}
        \centering
        \includegraphics[width=\linewidth]{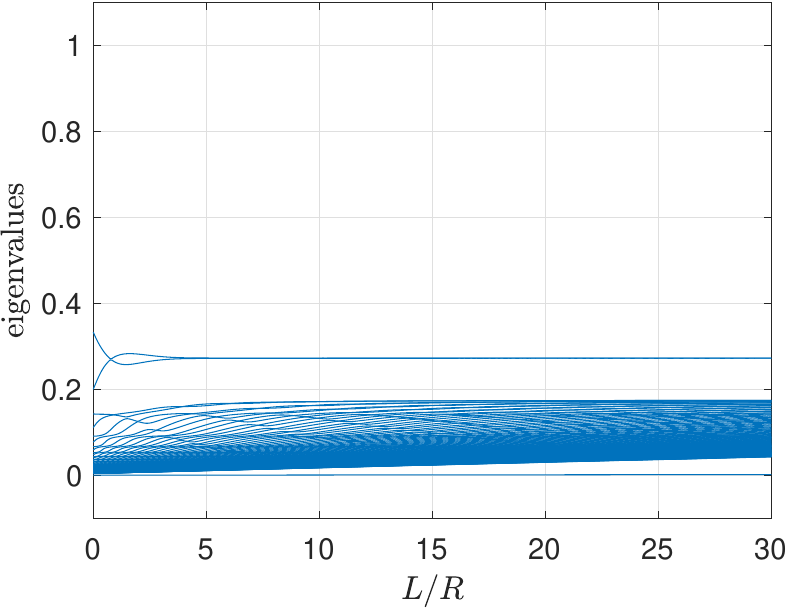}
        \caption{$m=1$}
        \label{fig:korv1}
    \end{subfigure}
    \caption{The sausage: computed eigenvalues as a function of $L/R$, where $R$
    is the radius of the semicircles and $L$ the length of the straight
    segments.}
    \label{fig:sausage}
\end{figure}

The sausage is obtained by rotating the cheese profile about a different axis. No
negative eigenvalues are detected; see Figure~\ref{fig:sausage}. The sausage has
no flat face, so Theorem~\ref{thm:flat} does not apply.

\begin{figure}[htbp]
    \centering
    \begin{subfigure}{0.48\textwidth}
        \centering
        \includegraphics[width=\linewidth]{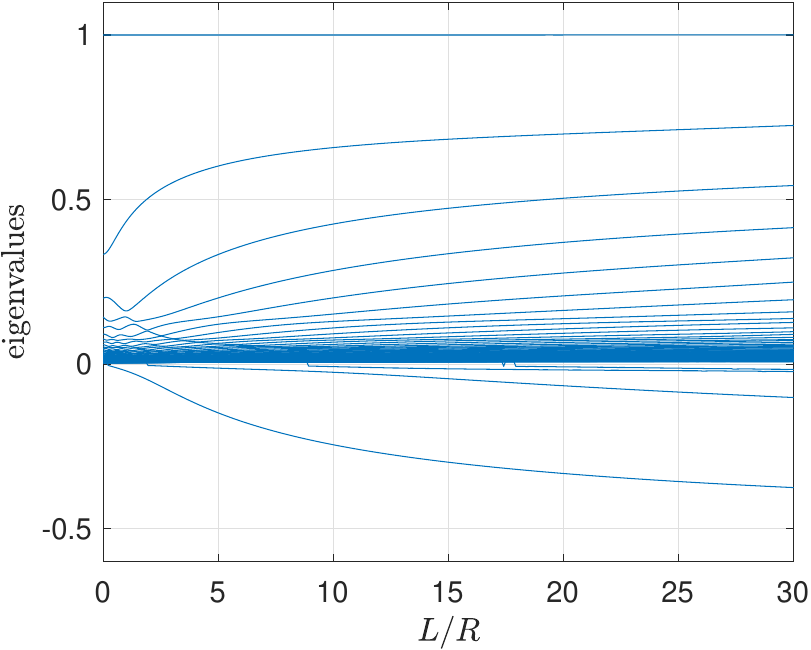}
        \caption{$m=0$}
        \label{fig:kokos0}
    \end{subfigure}
    \hfill
    \begin{subfigure}{0.48\textwidth}
        \centering
        \includegraphics[width=\linewidth]{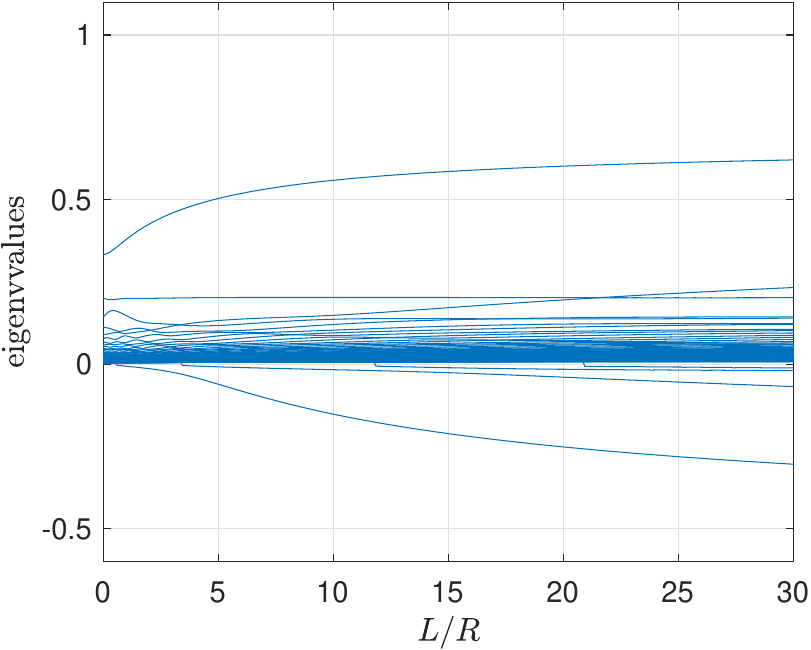}
        \caption{$m=1$}
        \label{fig:kokos1}
    \end{subfigure}
    \caption{The Swedish coconut macaroon: computed eigenvalues as a function of
    $L/R$, where $R$ is the radius of the arcs and $L$ the length of the straight
    segments.}
    \label{fig:coconut}
\end{figure}

The Swedish coconut macaroon has a flat face, so by Theorem~\ref{thm:flat}
infinitely many negative eigenvalues appear as soon as the length $L$ of the
straight segments that separate the arcs in the profile becomes positive; see
Figures~\ref{fig:curves}(\subref{fig:coconutcurve}) and~\ref{fig:coconut}.

\subsection{The superellipsoid}

\begin{figure}[htbp]
    \centering
    \begin{subfigure}{0.48\textwidth}
        \centering
        \includegraphics[width=\linewidth]{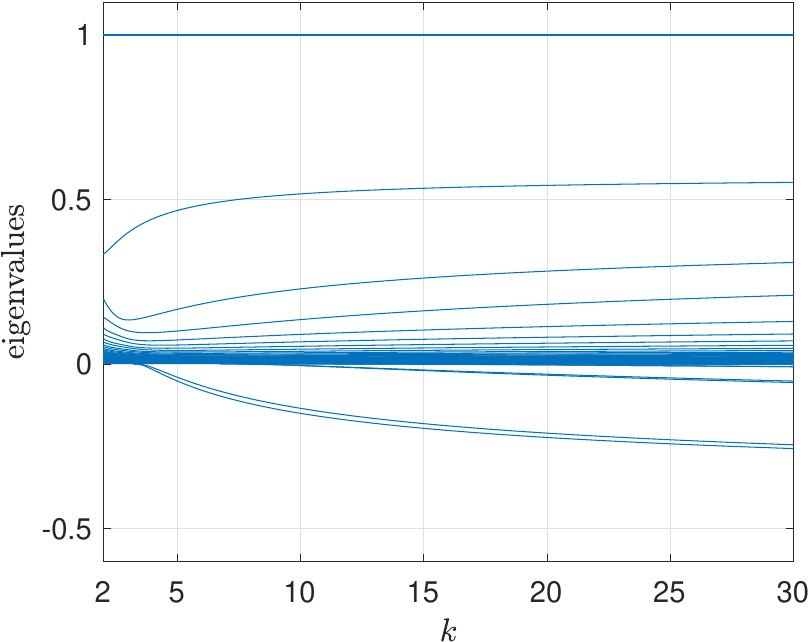}
        \caption{$m=0$}
        \label{fig:super0}
    \end{subfigure}
    \hfill
    \begin{subfigure}{0.48\textwidth}
        \centering
        \includegraphics[width=\linewidth]{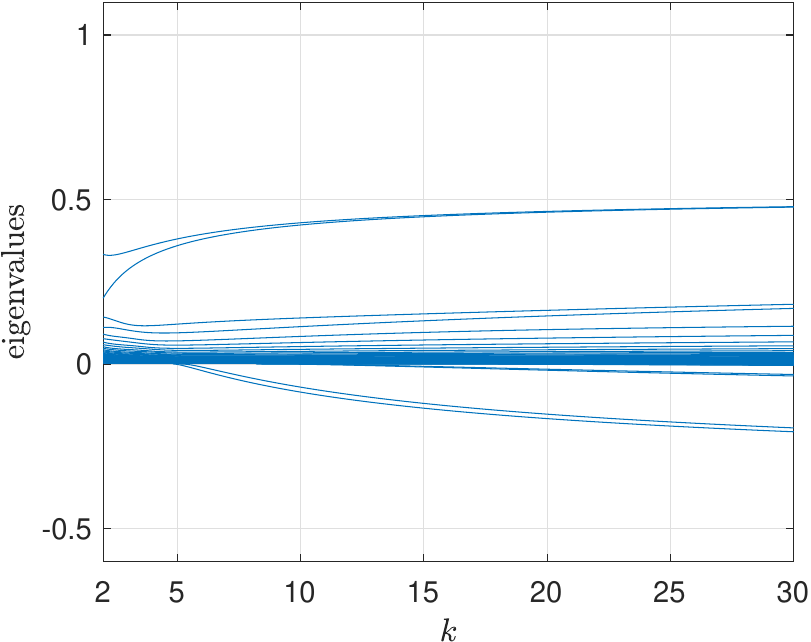}
        \caption{$m=1$}
        \label{fig:super1}
    \end{subfigure}
    \caption{Computed eigenvalues for the
    superellipsoid~\eqref{eq:superellipsoid} as a function of the exponent $k$.}
    \label{fig:super}
\end{figure}

For $k\geq2$ the superellipsoid is described by
\begin{equation}\label{eq:superellipsoid}
 (x_1^2+x_2^2)^{k/2}+\abs{x_3}^k=1.
\end{equation}
This is the surface of revolution obtained by rotating a superellipse about the
$x_3$-axis. At $k=2$ this is the unit sphere. As $k$ increases, the surface tends
to a cylinder. For $k>2$ both principal curvatures vanish at the poles and one
vanishes along the equator; for exponents~$k$ that are not even integers, the
regularity of the surface is of finite order. More and more negative eigenvalues
emerge as $k$ grows; see Figure~\ref{fig:super}. In
Example~\ref{ex:superellipsoid-cap} we showed that there exists at least one
negative eigenvalue if $k>4$, but as seen in Figure~\ref{fig:super} this is not
sharp.

\begin{figure}[htbp]
    \centering
    \includegraphics[height=5.75cm]{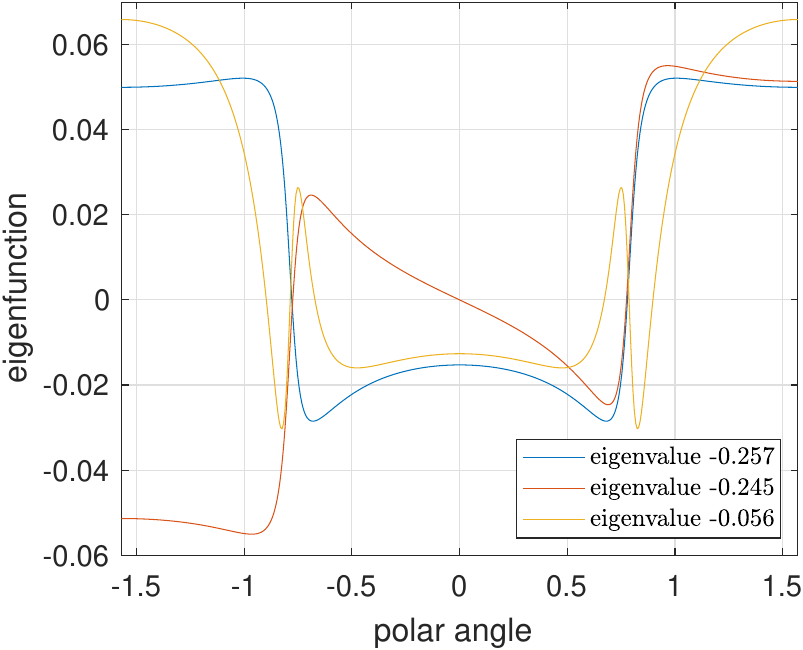}
    \caption{ Three eigenfunctions for the
    superellipsoid~\eqref{eq:superellipsoid} with $m=0$ and $k=30$, corresponding
    to negative eigenvalues, plotted as functions of the polar angle.}
    \label{fig:supereigf}
\end{figure}

For $m=0$ and $k=30$, Figure~\ref{fig:supereigf} shows three eigenfunctions
associated with negative eigenvalues. The eigenfunctions vary rapidly where the
curvature is large.

\clearpage

\appendix

\section{A warning about conventions}
\label{secwarning}

There is some variation in the literature concerning signs and normalizations of
layer potentials and of the NP operator; see, for
instance,~\cite[Section~2.1]{AKMP21}, where this lack of a uniform convention is
explicitly noted. Let $n$ denote the outward unit normal as above. One
convention, common in the literature on conductivity problems and plasmonic
resonance, is to let $\Gamma$ be the fundamental solution of $\Delta$, so that
$\Delta\Gamma=\delta_0$ (and hence $\Gamma(x)=-1/(4\pi\abs{x})$ in dimension three),
and to define (for $x\not\in \partial\Omega$)
\begin{equation}\label{singledouble}
 \mathcal S\varphi(x)=\int_{\partial\Omega}\Gamma(x-y)\varphi(y)\dd S_y,
 \qquad
 \mathcal D\varphi(x)=\int_{\partial\Omega}
 [\partial_{n_y}\Gamma(x-y)]\varphi(y)\dd S_y.
\end{equation}
Both the principal-value boundary operators $K$, $K^*$ and the operators
\[
 \tilde{K}\varphi(x)= \operatorname{p.v.}\!\int_{\partial\Omega}
 [\partial_{n_y}\Gamma(x-y)]\varphi(y)\dd S_y,
 \qquad
 \tilde{K}^*\varphi(x)=\operatorname{p.v.}\!\int_{\partial\Omega}
 [\partial_{n_x}\Gamma(x-y)]\varphi(y)\dd S_y
\]
are commonly referred to as NP operators. Note that $K=2\tilde{K}$. We call $K^*$, as
defined in~\eqref{eq:intro-np}, the NP operator; it is natural for transmission
problems since
\begin{equation}
\label{eq:jump}
  \partial_n(\mathcal S\varphi)\bigr|_{\pm}
   =\frac12(\pm I+K^*)\varphi.
\end{equation}
Here the $\pm$ indicates the limit of the values at $x+tn_x$ as $t\to 0^\pm$,
i.e., the $+$ indicates an exterior limit and the $-$ indicates an interior
limit. This is, for example, the convention used in~\cite{KangKimLeeShinYu2016}
and in much of the work of Ammari and Kang. In another widespread
boundary-integral convention one instead uses the fundamental solution
$\Phi=-\Gamma$ of $-\Delta$, which is positive in dimension three, and defines
the direct double-layer (or NP) operator by
\[
 D\varphi(x)=\operatorname{p.v.}\!\int_{\partial\Omega}
 [\partial_{n_y}\Phi(x-y)]\varphi(y)\dd S_y.
\]
With the same orientation of the normal this operator is $D=-\tilde{K}$ in the
preceding convention; see, for
example,~\cite{ChandlerWildeHaggerPerfektVirtanen2023}, where the competing sign
conventions are also explicitly discussed. Finally, a normalization going back to
the classical NP literature and used, for example, in~\cite{AKMP21,KPS07}
incorporates an additional factor $2$ and defines the double-layer operator by
\[
 K\varphi(x) = -2\,\operatorname{p.v.}\!\int_{\partial\Omega}
 [\partial_{n_y}\Phi(x-y)]\varphi(y)\dd S_y,
 \qquad \Delta \Phi=-\delta_0.
\]

Since $K=2\tilde{K}$ relative to the first convention above, the eigenvalues of
$\tilde{K}$ lie in $[-1/2,1/2]$, whereas those of $K$ lie in $[-1,1]$, and the
distinguished eigenvalue $1$ of $K^*$ becomes $1/2$ for $\tilde{K}^*$. Reversing
the choice of normal changes the sign of the double-layer and NP operators once
more. Consequently, when comparing spectral statements---in particular statements
concerning positive versus negative eigenvalues---it is important to specify the
sign of the fundamental solution, the orientation of the normal, the sign chosen
in the definition of the double layer, the possible factor $2$, and whether
\enquote {the NP operator} denotes the direct operator $K$ or its adjoint $K^*$.

\bibliographystyle{amsplain}
\bibliography{ourbib}

\end{document}